%% file: renorm_wp.tex
\documentclass{article}

\usepackage{amsmath}
\usepackage{amsfonts}
  \usepackage{enumerate}
\usepackage{mathptmx}

\usepackage{mathrsfs}
\usepackage{verbatim}
\usepackage[utf8]{inputenc}
\usepackage[T1]{fontenc}

\usepackage{ae,aecompl}
\usepackage{braket}
\usepackage{color}
\input{./macros}

\renewcommand{\hbar}{\bar{{\mathbb H}}^3}

\newcommand{\R}{\mathbb R}

\newcommand{\Rvol}{ V_R}
\newcommand{\Cvol}{V_C}
\newcommand{\T}{{\operatorname{Teich}}}

\def\eproof{$\Box$ \medskip}

\newcommand{\pslt}{\mathsf{PSL}_2(\mathbb C)}

\newcommand{\geod}{{\rm geod}}
\renewcommand{\cG}{{\mathcal G}}

\widemargins

\begin{document}

\title{\bf \Large 
  The Weil-Petersson gradient flow of renormalized volume and 3-dimensional convex cores} \author{Martin
   Bridgeman\thanks{M. Bridgeman's research  supported by NSF grants DMS-1500545, DMS-1564410}, \  Jeffrey
  Brock\thanks{J. Brock's research supported by NSF grants
    DMS-1608759, DMS-1849892}, \ and
  Kenneth Bromberg\thanks{K. Bromberg's research supported by NSF grants
    DMS-1509171, DMS-1906095}}

%\date{February 29, 2020}

\maketitle

\begin{abstract}
\input{abstract}
\end{abstract}

%\tableofcontents

\input{intro}

%\input{flow}
\input{background}

\input{completion}

\input{smallL2}

\input{bersslice}

\input{endgame}

%\input{noboundary}
\input{appendix}

\bibliography{bib,math}
\bibliographystyle{math}

\bigskip

{\sc
Boston College

\bigskip

Yale University

\bigskip

University of Utah

}
\end{document}

%% file: macros.tex
\newcommand{\bb}{\mathbb}

\usepackage{euscript}

\newcommand{\cx}{{\bb C}}

\newcommand{\widemargins}{
\setlength{\textwidth}{5.8in}
\setlength{\oddsidemargin}{0.25in}
\setlength{\evensidemargin}{0.25in}
}

\renewcommand{\bold}[1]{\medskip \noindent {\bf \boldmath #1
                        }\nopagebreak[4]}

\newcommand{\bdry}{\partial}

\newcommand{\del}{\partial}

\newcommand{\includesin}{\hookrightarrow}

\newcommand{\chat}{\widehat{\cx}}

\newcommand{\area}{\operatorname{area}}

\newcommand{\inj}{\operatorname{inj}}

\newcommand{\interior}{\operatorname{int}}
\renewcommand{\Im}{\operatorname{Im}}

\renewcommand{\Re}{\operatorname{Re}}

\newcommand{\sech}{\operatorname{sech}}

\newcommand{\Teich}{\operatorname{Teich}}

\newcommand{\vol}{\operatorname{vol}}

\newtheorem{theorem}{Theorem}[section]
\newtheorem{prop}[theorem]{Proposition}
\newtheorem{lemma}[theorem]{Lemma}
\newtheorem{cor}[theorem]{Corollary}
\newtheorem{conj}[theorem]{Conjecture}
\newtheorem{question}[theorem]{Question}

\newtheorem{bigtheorem}{Theorem}

\newcommand{\cC}{{\cal C}}

\newcommand{\cG}{{\cal G}}

\newcommand{\calC}{{\mathcal C}}

\newcommand{\calS}{{\mathcal S}}

%% file: abstract.tex
%% !TEX root =renorm_wp.tex

In this paper, we use the Weil-Petersson gradient flow for
renormalized volume to study the space $CC(N;S,X)$ of convex cocompact
hyperbolic structures on the relatively acylindrical 3-manifold
$(N;S)$. Among the cases of interest are the deformation space of an
acylindrical manifold and the Bers slice of quasi-Fuchsian space
associated to a fixed surface.  To
treat the possibility of degeneration along flow-lines to peripherally
cusped structures, we introduce a surgery procedure to yield a
surgered gradient flow that limits to the unique structure $M_\geod
\in CC(N;S,X)$ with totally geodesic convex core boundary facing
$S$. Analyzing the geometry of structures along a flow line, we show
that if $V_R(M)$ is the renormalized volume of $M$, then
$V_R(M)-V_R(M_\geod)$ is bounded below by a linear function of the
Weil-Petersson distance $d_{\rm WP}(\partial_c M, \partial_c
M_\geod)$, with constants depending only on the topology of $S$. The
surgered flow gives a unified approach to a number of problems in the
study of hyperbolic 3-manifolds, providing new proofs and
generalizations of well-known theorems such as Storm's result that
$M_\geod$ has minimal volume for $N$ acylindrical and the second
author's result comparing convex core volume and Weil-Petersson
distance for quasifuchsian manifolds.

%% file: intro.tex
%% !TEX root =renorm_wp.tex

\section{Introduction}

The use of a {\em geometric flow}, or a flow on a space of metrics on
a given manifold, has provided an abundantly fruitful approach to
understanding a manifold's structure. In our previous work
\cite{BBB}, we introduced a new geometric flow on the space of {\em
  hyperbolic} metrics on a 3-manifold that admits a hyperbolic
structure, showing how the flow can be used to discover the metric of
least {\em convex core volume}. In the present paper, we illustrate
how this flow provides an analytic version of results on convex core
volume that were available previously only through combinatorial
methods, demonstrating how this approach allows for conjectured
extensions to much more general cases.

When a hyperbolic 3-manifold $M$ admits a compact convex submanifold
we say it is {\em convex co-compact}; the geometry of the smallest
such submanifold, its {\em convex core}, carries all the interesting
information about its geometry. For such $M$ (or more generally
conformally compact Einstein manifolds), work of Graham and Witten
(\cite{GW}) in physics led to an alternative notion of {\em
  renormalized volume}. From a mathematical perspective, this concept
has been elaborated in a series of papers (see \cite{TT,ZT,KS08,KS12})
of Takhtajan, Zograf, Teo, Krasnov, and Schlenker. The renormalized
volume $\Rvol(M)$ of $M$ connects many analytic notions from the
deformation theory to the geometry of $M$ and is closely related to
classical objects such as the convex core volume $V_C(M)$ and the
Weil-Petersson geometry of Teichm\"uller space.

If $N$ is a compact 3-manifold admitting a complete hyperbolic
structure of finite volume, the renormalized volume gives an analytic
function $V_R:CC(N) \rightarrow \R$ where $CC(N)$ is the deformation
space of convex co-compact structures on $N$. We will give a precise
definition of $V_R$ later in the paper, but knowledge of its basic
properties will be largely sufficient for our purposes. In particular,
the differential $dV_R$ on $CC(N)$ is described in terms of the
classical {\em Schwarzian derivative} and can be used as a definition of $V_R$.

A convex co-compact structure $M \in CC(N)$ is naturally compactified
by a {\em complex projective structure} on $\partial N$. The
underlying conformal structure is the {\em conformal boundary}
$\partial_c M$ of $M$. The Schwarzian derivative associated to the
projective structure determines a holomorphic quadratic differential
$\phi_M \in Q(\del_c M)$. The utility of the renormalized volume
function lies in a particularly clean formula for its derivative,
first shown by Takhtajan-Zograf (\cite{ZT}) and Takhtajan-Teo
(\cite{TT}). A new proof was given by Krasnov-Schlenker (\cite[Lemma
8.5]{KS12}) using methods that are more closely aligned with the
present work. To state the result, we recall that $CC(N)$ is (locally)
parameterized by $\Teich(\partial N)$ and the cotangent space at
$\partial_c M$ is parameterized by $Q(\partial_c M)$. We then have:
\begin{theorem}[\cite{ZT,TT,KS12}]
\label{variational}
Let $\mu$ be an infinitesimal Beltrami differential on $\del_c M$. Then
$$d\Rvol(\mu) = \Re \int_{\del_c M} \phi_M\mu.$$
\end{theorem}

By integrating this formula along a Weil-Petersson geodesic and
applying the classical Kraus-Nehari bound on the $L^\infty$-norm of
$\phi_M$, Schlenker (\cite[Theorem 1.2]{schlenker-qfvolume}) obtained
the following for the quasifuchsian structure $Q(X,Y)$ on $N =
S\times[0,1]$ with conformal boundary $X\sqcup Y$:  
$$\Rvol(Q(X,Y)) \leq 3\sqrt{\frac{\pi}{2}|\chi(S)|}\, d_{\rm WP}(X,Y).$$
Furthermore, Schlenker showed that for quasifuchsian manifolds the
renormalized volume and the volume of the convex core are boundedly
related. A more refined version (see \cite[Theorems 2.16 and 3.7]{BBB}) is, 
$$\Cvol (Q(X,Y)) - 6\pi|\chi(S)| \leq \Rvol(Q(X,Y)) \leq \Cvol(Q(X,Y)).$$Combined, these gave a new proof of an upper bound on the volume the 
convex core of $Q(X,Y)$ in terms of $d_{\rm WP}(X,Y)$ originally due
to the second author \cite{Brock:wp}, resulting also in new approaches
to the study of volumes of fibered 3-manifolds 
in \cite{Brock:Bromberg:vol, Kojima:McShane:volumes} generalizing and sharpening known
estimates \cite{Brock:3ms1}.

Here, the variational formula (Theorem \ref{variational}) will be our
jumping off point to study the Weil-Petersson gradient flow of
$V_R$. It will be useful to restrict $V_R$ to certain subspaces of the
space of convex co-compact structures $CC(N)$. In particular, let
$(N;S)$ be a pair where $N$ is a compact hyperbolizable 3-manifold and
$S \subseteq \partial N$ is a collection of components of the
boundary. Then $CC(N;S,X) \subseteq CC(N)$ is the space of convex
co-compact hyperbolic structures on $N$ where the conformal boundary
on the complement of $S$ is the fixed conformal structure $X$. The
pair $(N;S)$ is {\em relatively incompressible} if the inclusion
$S\hookrightarrow N$ is $\pi_1$-injective and {\em relatively
  acylindrical} if there are no essential cylinders with boundary in
$S$. Note that the second condition implies the first.

In this paper our focus will be when $(N;S)$ is relatively
acylindrical. The cases of greatest interest are 1) when $S = \partial
N$ and $N$ itself is acylindrical and 2) when $N = S \times [0,1]$ and
$CC(N; S\times\{1\}, X)$ is a {\em Bers slice} of the space of
quasifuchsian structures. 
One important feature of relatively acylindrical pairs is that the
deformation space $CC(N;S,X)$ has a unique hyperbolic structure
$M_{\rm geod}$ where the components of the convex core facing $S$ are
totally geodesic. The main application of our study of the gradient
flow is the following: 
\begin{bigtheorem}\label{main bound}
Let $CC(N,S;X)$ be a relatively acylindrical deformation space. There exists $A(S)$, depending only on the topology of $S$, and a universal constant $\delta$ such that
$$A(S) \left(d_{\rm WP}(\partial_c M_{\rm geod}, \partial_c M) - \delta\right) \le V_R(M) -V_R(M_{\rm geod})$$
\end{bigtheorem}

For a Bers slice $CC(N; S\times\{1\}, X)$, we have $M_{\rm geod} =
Q(X,X)$ and both the convex core and renormalized volume of this
Fuchsian manifold are zero. Applying the above comparison between
renormalized volume and convex core volume we obtain: 
\begin{bigtheorem}
Let $S$ be a closed surface of genus $g\ge 2$. Then we have
$$A(S)(d_{\rm WP}(X,Y) - \delta) \le V_C(Q(X,Y)) \le
3\sqrt{\frac{\pi}2 |\chi(S)|}\, d_{\rm WP}(X,Y) + 6\pi|\chi(S)|.$$
\label{qfmain}
\end{bigtheorem}

Schlenker's argument in the quasifuchsian case also applies to
relatively acylindrical manifolds, so we have for any $M$ and $M'$ in
$CC(N;S,X)$ the following:
$$V_R(M) - V_R(M')\le 3\sqrt{\frac{\pi}2 |\chi(S)|}\, d_{\rm WP}(\partial_c M, \partial_c M').$$
If we let $M_{\rm geod} = M'$ then we get an upper bound on the expression in Theorem \ref{main bound}. The comparison between renormalized volume and convex core volume also extends to acylindrical manifolds (or any manifold with incompressible boundary). 

\begin{bigtheorem}\label{acylindrical bounds}
Let $N$ be a hyperbolizable, acylindrical 3-manifold. Then
\begin{eqnarray*}
A(\partial N) \left(d_{\rm WP}(\partial_c M_{\rm geod}, \partial_c M)-
  \delta\right) & \le & V_C(M) -V_C(M_{\rm geod})\\
&\le& 3\sqrt{\frac{\pi}2 |\chi(\partial N)|}\, d_{\rm WP}(\partial_c
M_{\rm geod}, \partial_c M) + 3\pi|\chi(\partial N)|
\end{eqnarray*}
where $A$ and $\delta$ are  as in Theorem \ref{main bound}.
\label{maincorevolume}\end{bigtheorem}

\bold{Remark.} We note that the constants in Theorem~\ref{acylindrical
  bounds} {\em depend only on the topology of $\partial N$}. While we
expect the second author's original method combined with Thurston's
compactness theorem for hyperbolic structures on acylindrical
manifolds should also produce a similar bound, the constants in such
an approach would depend on the topology of $N$, due to the
application of Thurston's result. The approach taken here is thus not
only more direct but produces a stronger result. In particular,
while Thurston's compactness theorem implies that the convex core of
$M_{\rm  geod}$ has a bi-Lipschitz embedding into any complete
hyperbolic 
structure on $N$ where the bi-Lipschitz constants only depend on
$N$, it is natural to conjecture that these bi-Lipschitz constants
only depend on $\partial N$. Theorem \ref{acylindrical bounds} can be
taken as some evidence for this conjecture.

We note that a positive resolution of this conjecture would also imply
Minsky's conjecture that the diameter of the {\em skinning map} is
bounded by constants only depending on $\partial N$, and provide an
approach to improving related estimates for the models of
\cite{BMNS:bounded:models}.

\subsection{The Weil-Petersson gradient flow of renormalized volume}
One of the main purposes of this paper is to develop the structure theory of the gradient flow $V$ for renormalized volume $V_R$. 
From this development, the above results will follow directly. We show
that flow provides a powerful new tool to investigate the internal
geometry of ends of hyperbolic 3-manifolds.  

To give a basic outline of the main ideas of the paper, we begin with
a general discussion of gradient flows which we will then apply to the
gradient of renormalized volume. Let $f$ be a smooth function on a
non-compact Riemannian manifold $X$ and assume 
\begin{enumerate}[(a)]
\item $f$ is bounded below;\label{below}

\item the gradient flow of $f$ is defined for all time;

\item $\|\nabla f\|\le C$;\label{up}

\item $f$ has a unique critical point $\bar x$;\label{critical point}

\item for all $\epsilon>0$ there exists a $A>0$ such that if $d(x, \bar x) \ge \epsilon$ then $\|\nabla f\| \ge A$.
\label{small ball}
\end{enumerate}
By integrating $\|\nabla f\|$ along a distance minimizing path between points $x$ and $x'$ we immediately see that \eqref{up} implies that
$$|f(x) - f(x')| \le Cd(x,x').$$
Clearly, we cannot expect a similar lower bound to hold as the level
sets of $f$ may have infinite diameter. Instead, we obtain lower
bounds when $x'=\bar x$, the unique critical point. In particular, let
$x_t$ be a flow line of $-\nabla f$ with $x= x_0$. We then have 
$$f(x) - f(x_a)= \int_0^a \|\nabla f(x_t)\|^2 dt.$$
By \eqref{below}, $\underset{a\to \infty}{\lim} f(x_a)$ exists so as $a\to\infty$, the improper integral is convergent. Therefore there will be an increasing sequence of $t_i$ with $\|\nabla f(t_i)\|\to 0$ so, by \eqref{small ball}, the flow line $x_t$ will accumulate on $\bar x$. Fix some $\epsilon>0$ with corresponding $A>0$ as in \eqref{small ball} and let $I_\epsilon \subset [0,\infty)$ be those values $t$ where $d(x_t, \bar x) > \epsilon$. Then for $t\in I_\epsilon$ we have $\|\nabla f(x_t)\| \ge A$ and the length of the path $x_t$ restricted to $I_\epsilon$ will be at least $d(x,\bar x) - \epsilon$. Therefore
\begin{eqnarray*}
f(x) - f(\bar x) & =& \int_0^\infty \|\nabla f(x_t)\|^2 dt \\
& \ge & \int_{I_\epsilon} \|\nabla f(x_t)\|^2 dt \\
& \ge & A\int_{I_\epsilon} \|\nabla f(x_t)\| dt\\
& \ge & A(d(x, \bar x) - \epsilon)
\end{eqnarray*}
which gives the desired linear lower bound.

Unfortunately, when we replace $f$ with the renormalized volume function $V_R$, property \eqref{small ball} will not hold (but the others will). To mimic what happens in our generic setting, we let $\bar X$ be the metric completion of our Riemannian manifold $X$ and $\mathcal G \subset \bar X$ a subset. We replace \eqref{small ball} with the following three properties:
\begin{enumerate}[(\ref{small ball}-1)]
\item for all $\epsilon>0$ there exists a $A>0$ such that if $d(x, \mathcal G) \ge \epsilon$ then $\|\nabla f(x)\| \ge A$.
\label{small nbd}

\item there exists an $N>0$ such that in any subset of $\mathcal G$ with more than $N$ elements there are at least two that a distance $\delta_0$ apart;\label{separation}

\item for every $x_0 \in \mathcal G$ there is a path $x_t$ starting at $x_0$ with $x_t \in X$ for $t>0$ and $f(x_t)< f(x_0)$.
\label{decrease path}
\end{enumerate}
While the overall structure of the argument will remain the same, some modifications are necessary. First, we need to construct a {\em surgered flow} $x_t$ where
\begin{itemize}
\item $x_0 = x$;

\item the function $t\mapsto f(x_t)$ satisfies $f(x_t) < f(x_0)$;

\item outside of the $\epsilon$-neighborhood of $\mathcal G$, $x_t$ is the gradient flow;

\item $x_t \to \bar x$ as $t\to\infty$.
\end{itemize}
To construct $x_t$ we start the gradient flow at $x$. If it limits to $\bar x$ (as we conjecture it will for renormalized volume) then we are done. If not, we limit to some other point in $\mathcal G$. We reparameterize so that this happens in finite time and then use (\ref{small ball}-\ref{decrease path}) to restart the flow. If this converges to $\bar x$ we stop; if not we repeat. The first three bullets follow directly from this construction.

As before we fix an $\epsilon$ and $A$ as in (\ref{small ball}-\ref{small nbd}) and let $I_\epsilon(a) \subset [0,a]$ be those $t \in [0,a]$ where $d(x_t, \mathcal G) > \epsilon$. If $L_\epsilon(a)$ is the length of the path $x_{[0,a]}$ restricted to $I_\epsilon(a)$ then the above argument gives
$$f(x) -f(x_a) \ge A L_\epsilon(a).$$
A simple geometric argument, using (\ref{small ball}-\ref{separation}), shows that $L_\epsilon(a)$ grows linearly in both the number of points of $\mathcal G$ that $x_t$ passes through and in the distance $d(x, x_a)$. In particular, if $x_t$ passes through infinitely many points in $\mathcal G$ then $L_\epsilon(a) \to \infty$ as $a\to\infty$ so $f(a) \to -\infty$, contradicting \eqref{below}. Therefore $x_t$ only passes through finitely many points in $\mathcal G$ which implies that the surgered flow converges to the critical point. Therefore if we take the limit of the above inequality we have
$$f(x) - f(\bar x) \ge A L_\epsilon(\infty)$$
and as $L_\epsilon(\infty)$ is bounded below by a linear function of $d(x, \bar x)$ we have our bound.

We now apply this discussion to the renormalized volume function $V_R$ on a relatively acylindrical deformation space $CC(N;S,X)$. Properties (\ref{below})-(\ref{critical point}) are already known so we will focus on (\ref{small ball}-\ref{small nbd})-(\ref{small ball}-\ref{decrease path}). In particular, we need to understand when $\|\nabla V_R\|$ is small. By Theorem \ref{variational} we have that the Weil-Petersson gradient of $V_R$ is given by the {\em harmonic} Beltrami differential
$$\nabla V_R (M) = \frac{\overline{\phi_M}}{\rho_M}$$
where $\rho_M$ is the area form for the hyperbolic metric on $\partial_c M$ and $\phi_M$ is the quadratic differential associated to the projective structure  on  the components of $\partial_c M$ corresponding to $S$. The norm of $\nabla V_R$ is then the $L^2$-norm of $\phi_M$. This $L^2$-norm is zero exactly when $\phi_M = 0$. As $\phi_M$ is  the Schwarzian derivative of the univalent map uniformizing the components of $\partial_c M$ corresponding to $S$ (see \cite{KS12}), $\phi_M = 0$ implies that the uniformizing maps are M\"obius. It follows that if the norm of $\nabla V_R$ is zero then the components of the boundary of the convex core facing $S$ are totally geodesic. In a relatively acylindrical deformation space there is exactly one such manifold (which is why \eqref{critical point} holds) and one might hope that when $\|\phi_M\|_2$ is small we are near this critical point. If this were so \eqref{small ball} would hold. Unfortunately, it does not. While $\|\phi_M\|_2$ being small will imply that $M$ is near a hyperbolic manifold whose convex core boundary (facing $S$) is totally geodesic, this manifold may have {\em rank one cusps}. 

To state this more precisely, if $GF(N; S,X)$ is the space of {\em geometrically finite} hyperbolic structures on $(N;S,X)$,  then the map $M\mapsto \partial_c M$ is a bijection from $GF(N;S,X)$ to the Weil-Petersson metric completion $\overline{\Teich(S)}$ of Teichm\"uller space where points in the completion are {\em noded hyperbolic structures} on $S$ (see \cite{Masur:WP}). Nodes in the conformal boundary correspond to rank one cusps in the hyperbolic 3-manifold. The triple $(N;S,X)$ determines a subset $\mathcal G(N;S,X)$ of $\overline{\Teich(S)}$ where the corresponding hyperbolic structures have totally geodesic boundary facing $S$. With $\mathcal G = \mathcal G(N;S,X)$ defined we can briefly describe how we will verify (\ref{small ball}-\ref{small nbd})-(\ref{small ball}-\ref{decrease path}).

Property (\ref{small ball}-\ref{small nbd}) is the following theorem and its proof will occupy much of the paper:
\begin{bigtheorem}\label{small-L2}
For all $\epsilon>0$ there exists $A = A(\epsilon, S)$ such that if $M \in CC(N;S,X)$ with $\|\phi_M\|_2 \le A$ then there is an $M' \in \mathcal G(N;S,X) \subset GF(N;S,X)$ such that
$$d_{\rm WP}(\partial_c M, \partial_c M') \le \epsilon.$$
\end{bigtheorem}

Property (\ref{small ball}-\ref{separation}) follows from Wolpert's strata separation theorem (Theorem \ref{strata}). For a noded surface $Y \in \partial\overline{\Teich(S)}$, we denote the family of curves given by the nodes by $\tau_Y$. Then Wolpert's strata separation theorem implies  there is a universal constant $\delta_0> 0$ such that if $Y_1, Y_2 \in  \partial\overline{\Teich(S)}$ with geometric intersection $i(\tau_{Y_1},\tau_{Y_2}) \neq 0$, then $d_{\rm WP}(Y_1, Y_2) > \delta_0$. Thus (\ref{small ball}-\ref{separation})  holds with $N = 2^{n(S)}$ where $n(S)$ is the maximal number of disjoint simple closed curves on $S$ as any collection of greater than $N$ noded surfaces in $\partial\overline{\Teich(S)}$  contains two  that have intersecting nodes.

Finally Property (\ref{small ball}-\ref{decrease path}) follows by unbending the nodes by decreasing the bending angle from $\pi$ along the nodes to some angle $\theta < \pi$. Such a deformation was constructed by Bonahon-Otal (\cite{BO:bending}).
Using the variational formula for $V_R$ it can be easily shown that $V_R$ satisfies  Property (\ref{small ball}-\ref{decrease path}) along this path (see Proposition \ref{no_local_min}) as required.

\subsection{Constants}
A striking feature of Schlenker's proof of the second author's upper
bounds for volume is that the constants are very
explicit. Unfortunately we lack the same control of constants in
our lower bounds as there is one place in the proof, the use of
McMullen's contraction theorem for the skinning map, that we fail to
control constants explicitly. If we assume, optimistically, that the
contraction constant does not depend on the manifold then we can at
least understand the asymptotics. With this assumption the
multiplicative constant in our lower bound will decay exponentially
with exponent of order $g^2$, where $g$ is the genus. On the other
hand, the additive constants will decay to zero even without
controlling the contraction constant. This should be compared to work
of Aougab-Taylor-Webb (\cite{ATW:effective}) who produced an effective
lower bound in the quasifuchsian case via the second author's
combinatorial 
methods. Their multiplicative constants decay exponentially with
exponent of order $g \log g$ which is better than ours but their
additive constant grows, also of order $g\log g$, rather than decays.
\subsection{Questions and Conjectures} 
A central feature of the surgered gradient flow of $-V_R$ on a relatively acylindrical deformation space is that it converges to the unique structure whose convex core has totally geodesic boundary.
While in this paper we will focus on relatively acylindrical deformation spaces, the gradient flow is defined on the deformation space of any hyperbolizable 3-manifold as is a surgered flow. We conjecture:

\begin{conj}
The surgered gradient flow either converges to a hyperbolic structure whose convex core has totally geodesic boundary or it finds an obstruction to the existence of such a structure. More concretely, either
\begin{itemize}
\item $N$ is acylindrical and $M_t \to M_{\rm geod}$ or
\item there is an essential annulus or compressible disk whose boundary has small length in $\partial_c M_t$ for some $t$.
\end{itemize}
\end{conj}

In fact we expect that the surgeries are unnecessary. Here is a more concrete conjecture when the manifold has incompressible boundary.

\begin{conj}\label{window_frame}
  Let $N$ have incompressible boundary. Then for $M \in
  CC(N)$ the renormalized volume gradient flow $M_t$ starting at $M$
  has the property that for any simple closed curve $\gamma$ on $\bdry
  N$ the geodesic length $\ell_{M_t}(\gamma^*)$ tends to zero if and
  only if $\gamma$ lies in the window frame.
\end{conj}
See Thurston's paper \cite{Thurston:hype3} for the definition of the window of a hyperbolic 3-manifold with incompressible boundary. 

In effect, the renormalized volume gradient flow realizes the
geometric decomposition of the manifold into pieces by pinching
cylinders corresponding to the window boundary, cutting the convex
core of the manifold
into pared acylindrical pieces with totally geodesic boundary and
Fuchsian ``windows.''

Other questions relate to the internal geometric structure of convex
cocompact ends and how the flow relates to their internal
structure. To avoid technicalities, for the remainder of this section we will assume that our manifolds are acylindrical. 

Let $\mathcal C(M, L)$ the collection of simple closed curves on
$\partial M$ that have geodesic length $\le L$ in $M$ and let
$\mathcal F(M, L)$ the collection of simple closed curves on $\partial
M$ that have length $\le L$ on some $\partial_c M_t$ where $M_t$ is
the gradient flow starting at $M$.

\begin{question}
Given $L>0$ does there exist an $L'>0$ such that
$$\mathcal F(M, L') \subset \mathcal C(M,L)$$
and
$$\mathcal C(M, L') \subset \mathcal F(M,L) ?$$
\end{question}
A stronger version of this question is the following.
\begin{question} Does the flow give a continuous family of
  bi-Lipschitz 
  embeddings into the initial manifold? In other words, for $s <t$
does the convex core of $M_t$ embed in the convex core of $M_s$ in a
bi-Lipschitz manner?
\label{question:bilip:embedding}
\end{question}
Note that a positive answer to this question would have
applications. First, it would imply Thurston's compactness theorem for
deformation spaces of acylindrical manifolds. A suitable
generalization of this conjecture to the general incompressible case
would also imply Thurston's relative compactness theorem in this
setting. It would also imply the following conjecture that was
mentioned above: 
\begin{conj}
Let $N$ be an acylindrical 3-manifold. Then for all $M \in CC(N)$ the
convex core of $M_{\rm geod}$ has a bi-Lipschitz embedding in $M$ with
constants only depending on $\partial N$. 
\end{conj}
We note that as gradient flow lines are Weil-Petersson
quasi-geodesics, relative stability properties established in
\cite{Brock:Masur:pants} for low-genus cases (genus two or lower
complexity) for such quasi-geodesics would control the behavior of
manifolds along the flow $M_t$ when $\bdry N$ has genus two. This
observation gives an approach to
Question~\ref{question:bilip:embedding} in such cases. Such stability
fails to hold in higher genus cases, so other properties of the flow
would be required. The question is reminiscent of similar questions
involving the relation of Weil-Petersson geodesics to properties of
ends of hyperbolic 3-manifolds and the models of
\cite{Brock:Canary:Minsky:elc}.

\bold{Acknowledgement.} We would like to thank MSRI and Yale
University for their hospitality while portions of this work were
being completed. We also thank Dick Canary, Curt McMullen, and Yair
Minsky for helpful conversations.

%% file: background.tex
%% !TEX root =renorm_wp.tex

\section{Background and notation}
In what follows, we fix $S$ to be a closed orientable surface with connected components having  genus at least two.
\subsubsection*{Norms on quadratic differential and metrics on Teichm\"uller space}
Let $\Omega^{p,q}(Y)$ be the space of $(p,q)$-differentials on a Riemann surface $Y$. Given $\phi\in\Omega^{2,0}(Y)$ (a {\em quadratic differential})  and $\mu\in\Omega^{-1,1}(Y)$ (a {\em Beltrami differential}) the product $\mu\phi$ is $(1,1)$-differential which can canonically be identified with a $2$-form so we have a pairing
$$\langle \phi,\mu\rangle = \int_Y \mu\phi.$$
In particular these two spaces are naturally dual.

We also have the subspace  $Q(Y) \subset \Omega^{2,0}(Y)$ of {\em holomorphic quadratic differentials}. This space is important as it is canonically identified with the co-tangent space $T^*_Y\Teich(S)$. The tangent space $T_Y\Teich(S)$ is then a quotient of $\Omega^{-1,1}(Y)$. In particular define 
$$N(Y) = \{\mu\in\Omega^{-1,1}(Y) | \langle \phi,\mu\rangle  = 0 \mbox{ for all } \phi \in Q(Y)\}$$
and then
$$T_Y\Teich(S) = \Omega^{-1,1}(Y)/N(Y).$$

If $\rho_Y$ is the area form for the hyperbolic metric on $Y$ and $\phi\in\Omega^{2,0}(Y)$ then $|\phi|/\rho_Y$ is also a function and we define $\|\phi(z)\| = |\phi(z)|/\rho_Y(z)$ to be the pointwise norm. We let $\|\phi\|_p$ be the $L^p$-norm of this function on $Y$, again with respect to the hyperbolic area form. Given $\mu \in \Omega^{-1,1}(Y)$ we define the $L^q$-norm (with $1/p+1/q=1$) of the equivalence class $[\mu]\in T_Y\Teich(S)$ by
$$\|[\mu]\|_q = \underset{\phi \in Q(Y)\backslash\{0\}}{\sup} \frac{\left|\langle \phi, \mu\rangle\right|}{\|\phi\|_p} \le \|\mu\|_q.$$
For $p=1$ this norm on $T_Y\Teich(S)$ gives the {\em Teichm\"uller metric} on $\Teich(S)$ and for  $p=2$ it gives the Weil-Petersson metric. Note that the Teichm\"uller metric is a Finsler metric while the Weil-Petersson metric is Riemannian as the $L^2$-norm on $Q(Y)$ can be given as an inner product. In particular, the $L^2$-norm on $Q(Y)$ is given by the inner product
$$(\psi, \phi) = \Re \int_Y \psi\bar\phi/\rho_Y.$$
From this we see that if $f\colon \Teich(S) \to \R$ is a smooth function then its differential $df$ is an assignment of a holomorphic quadratic differential $\phi_Y$ to each $Y\in \Teich(S)$. Its Weil-Petersson gradient is the vector field is represented  at each $Y$ by a Beltrami differential $\mu_Y$ where for all $\psi \in Q(Y)$ we have
$$(\psi, \phi_Y) = \langle \psi, \mu_Y \rangle.$$
It is a standard fact (and not hard to check directly) that $[\mu_Y]$ is represented by the {\em harmonic Beltrami differential} $\bar\phi_Y/\rho_Y$ and that
$$\|[\mu_Y]\|_2 = \|\bar\phi_Y/\rho_Y\|_2 = \|\phi_Y\|_2.$$

\subsubsection*{Collars}
We  state the Collar lemma originally due to Keen (\cite{Keen:collar}). We give it in a  form due to Buser \cite{Buser:collars}.

\begin{theorem}[Buser, {\cite{Buser:collars}}]
Let $Y$ be a  complete hyperbolic surface and $\gamma$ a simple closed geodesic of length $\ell_\gamma(Y)$. Then the collar $B(\gamma)$ of width $w(\gamma) =  \sinh^{-1}\left(1/\sinh\left(\frac{\ell_{\gamma}(Y)}{2}\right)\right)$ is embedded. If $z \in B(\gamma)$, then
$$\sinh\left( {\rm inj}_Y(z)\right) = \sinh\left(\ell_{\gamma}(Y)/2\right) \cosh \left( d(z,\gamma) \right).$$
Furthermore for any two disjoint geodesics the  collars are disjoint.
\end{theorem}

If $\ell_\gamma(Y) \leq 2\epsilon_2$ then we define the {\em standard collar} of $\gamma$ as 
$$\{ z \in B(\gamma)\ | \ {\rm inj}_Y(z) \leq \epsilon_2\}.$$
We note that it follows from  the collar lemma (see \cite{Buser:collars}) that the standard collar consists of all points in $Y$ that lie on a curve of length  $\leq 2\epsilon_2$ which is homotopic to $\gamma$.
 
For $S$ a finite type surface we  define $n(S)$ to be the maximal number of disjoint simple closed curves in $S$. For $S$ a  surface of genus $g$ and $k$ punctures we have $n(S) = 3g-3+k$ and for $S$ with connected components $S_i$ then $n(S) = \sum_i n(S_i)$.

\subsubsection*{Hyperbolic 3-manifolds}
Let $(N,P)$ be a pared 3-manifold (see
e.g. \cite{thurston-geometrization}) and $S$ a collection of
components 
of $\partial N - P$. Then the triple $(N,P;S)$ is {\em relatively
  acylindrical} if no essential cylinder has boundary 
 in S.
The acylindricity condition implies that all components of $S$ are incompressible.

A complete hyperbolic 3-manifold $M$ on the interior of $N$ naturally has the structure of a pared 3-manifold. This is simplest to describe when $M$ is geometrically finite and, as this is the only setting we will consider, we stick to this case. Let $\bar M$ be the union of $M$ and its conformal boundary. Then there is a paring locus $P \subset \partial N$ such that $\bar M$ is homeomorphic to $N - P$. The paring locus $P$ is a collection of annuli and tori. These are the {\em rank one} and {\em rank two cusps} of $M$. In particular, a curve in $M\subset N$ has parabolic holonomy if and only if it is homotopic into $P$.

Let $MP(N,P)$ be the space of geometrically finite hyperbolic structures on the interior of $N$ with induced pared manifold structure $(N,P)$. (These are {\em minimally parabolic} structures on $(N,P)$ - every parabolic is contained in $P$.) Now fix a conformal structure $X$ on the complement of $S$ in
$\partial N - P$ and let $MP(N,P;S,X)\subset MP(N,P)$ be those hyperbolic structures with
conformal boundary $X$ on the complement of $S$. Then by the deformation
theory of Kleinian groups (see e.g. \cite{Kra:crash}) we have the
parametrization $MP(N,P;S,X) 
\simeq \T(S)$.
The space $MP(N,P;S,X)$ is a {\em quasi-conformal deformation space}; any two hyperbolic manifolds in $MP(N, P; S,X)$ are quasi-conformal deformations of each other with the deformation supported on $S$.

Our results on renormalized volume will only apply to manifolds where $P$ is empty. However, in the course of the proof it will be necessary to consider hyperbolic 3-manifolds with cusps.

\subsubsection*{Schwarzian derivatives and projective structures}
Let $f\colon\Delta\to \chat$ be a locally univalent map on the unit disk $\Delta \subset \cx$. The Schwarzian derivative is the quadratic differential given by
$$Sf(z) = \left(\left(\frac{f''(z)}{f'(z)}\right)'- \frac12\left(\frac{f''(z)}{f'(z)}\right)^2\right)dz^2.$$
If $f$ is a M\"obius transformation then $Sf = 0$, and in general, $Sf$ measures how much $f$ differs from a M\"obius transformation. We also have the following composition rule
$$S(f\circ g) (z) = Sf(g(z))g'(z)^2 + Sg(z).$$
Observe that if $f$ is a M\"obius transformation then $S(f\circ g)  = Sg$, while if $g$ is a M\"obius transformation $S(f\circ g)(z)  = Sf(g(z))g'(z)^2$.

Let $\Gamma$ be a Fuchsian group such that $Y = \Delta/\Gamma$. A {\em projective structure} on $Y$ is given by a locally univalent map $f\colon\Delta\to\chat$ (the {\em developing map}) with a {\em holonomy representation} $\rho\colon \Gamma\to \pslt$ such that for all $\gamma \in \Gamma$ we have
$$f\circ \gamma = \rho(\gamma) \circ f.$$
The composition rule for the Schwarzian implies that $Sf$ descends to a a holomorphic quadratic differential in $Q(Y)$.

%% file: completion.tex
%% !TEX root =renorm_wp.tex

\subsubsection*{The Weil-Petersson completion and its stratification}

While the Teichm\"uller metric is complete, there are paths with
finite length in the Weil-Petersson metric that leave every compact
subset of Teichm\"uller space. Our goal in this section is to describe
some of the basic structure of the {\em completion} of the
Weil-Petersson metric. Points in this metric completion are naturally
parametrized by families of {\em Riemann surface with nodes}, namely,
a degeneration of a finite area hyperbolic Riemann surface obtained
by collapsing the curves in a multicurve to cusps.

Given a compact surface $S$,  the {\em complex of curves}
$\calC(S)$ is a simplicial complex organizing the isotopy classes of
simple closed curves on $S$ that do not represent boundary components.
To each isotopy class $\gamma$ we associate a vertex $v_\gamma$, and
each $k$-simplex $\sigma$ is the span of $k+1$ vertices whose
associated isotopy classes can be realized disjointly on $S$.

It is due to Masur \cite{Masur:WP} that the completion of $\Teich(S)$
with the Weil-Petersson metric is identified with the {\em augmented
  Teichm\"uller space}, obtained by adjoining at infinity the Riemann surfaces
with nodes. A point in the completion is given by a choice
of the multicurve $\tau$, a ($0$-skeleton of a) simplex in $\calC(S)$, and finite area hyperbolic structures on the
complementary subsurfaces $S \setminus \tau$. The completion is
stratified by the simplices of $\calC(S)$: the collection of noded
Riemann surfaces with nodes determined by a given simplex $\sigma$
lies in a product of lower-dimensional Teichm\"uller spaces determined
by varying the structures on $S \setminus \tau$. This {\em
  stratum} of the completion, $\calS_\tau$, inherits a natural metric
from the Weil-Petersson metric, which by Masur (see \cite{Masur:WP})  is isometric to
the product of Weil-Petersson metrics on the Teichm\"uller spaces of
the complementary subsurfaces.

The Teichm\"uller space, with this `augmentation' by its
Weil-Petersson completion, naturally descends under the action of the
mapping class group to a finite diameter metric on the Deligne-Mumford
compactification of the moduli space of Riemann surfaces. If $\overline\Teich(S)$ is the completion 
then we can describe the strata as follows
$$\calS_\tau = \{ X \in \overline\Teich(S)\ |\  \ell_\gamma(X) =0 \mbox{ if and only if } \gamma \in \sigma\}$$
where $\ell_\gamma$ is the extended length function of $\gamma$. 

We note
that if $\tau_0 \subseteq \tau_1$ are simplices in $\calC(S)$, then
we have $\calS_{\tau_1} \subseteq \overline{\calS_{\tau_0}}$. 

In his investigation of the geometry of the completion, Wolpert showed  the following; 
\begin{theorem}[Wolpert, {\cite[Cor. 22]{Wolpert:compl}}]\label{node_bound}
There is a positive constant $\delta_0$ so that  either $i(\tau_0,\tau_1) = 0$ and the closures of the strata
$\calS_{\tau_0}$ and $\calS_{\tau_1}$ intersect or $i(\tau_0,\tau_1) > 0$
$$d_{\rm WP}(\calS_{\tau_0},\calS_{\tau_1}) \ge \delta_0.$$
\label{strata}
\end{theorem}

We note that a simple calculation shows that $\delta_0 > 6$ (see  \cite{BB:strat}).

\section{Hyperbolic 3-manifolds with small Schwarzian derivative}
Before proving Theorem \ref{small-L2} we set some notation.
Let $(N,P;S)$ be a relatively acylindrical triple where $P$ is a collection of tori and $X$ a conformal structure on the complement of $S$ in $\partial N - P$.  We consider the following;

\begin{itemize}
\item $\tau$ is a simplex in $\mathcal C(S)$.

\item $P_\tau$ is the union of $P$ and the curves in $\tau$.

\item $S_\tau$ is the complement of $\tau$ in $S$.
\end{itemize}
Note that the new triple $(N, P_\tau; S_\tau)$ is still relatively acylindrical and the complement of $S_\tau$ in $\partial N - P_\tau$ is homeomorphic to the complement of $S$ in $\partial N - P$. We then have
$$GF(N,P;S,X) = \underset{\tau}{\sqcup} MP(N, P_\tau; S_\tau, X)$$
so $GF(N,P;S,X)$ is naturally parameterized by the Weil-Petersson completion $\overline{\Teich(S)}$ of Teichm\"uller space.

We next set:
\begin{itemize} 
\item If $Y \in \overline{\Teich(S)}$ then $M_Y$  is the hyperbolic manifold in $GF(N,P;S,X)$ under the above identification $GF(N,P;S,X) \cong \overline{\Teich(S)}$.

\item $\phi_Y$ is the Schwarzian quadratic differential given by the projective structure on $Y$ induced by $M_Y$.

\end{itemize}
We are especially interested in those manifolds in $GF(N,P; S, X)$ where the boundary of the convex core facing $S$ is totally geodesic. We set notation for this set:
\begin{itemize}
\item $Y_\geod^\tau$ is the unique conformal structure in $\Teich(S_\tau)$ such that the component of the convex core of $M_{Y^\tau_\geod}$ facing $S_\tau$ is totally geodesic.

\item $\cG(N, P;S,X)$ is the union of the $Y_\geod^\tau$.

\item If $\tau = \emptyset$ then we set $Y_\geod = Y^\tau_\geod$ and $M_\geod = M_{Y_\geod}$.
\end{itemize}

We have the following elementary observation.

\begin{lemma}
Let $(N,P;S)$ be a relatively acylindrical triple where $P$ is a collection of tori and $X$ a conformal structure on the complement of $S$ in $\partial N - P$. Then the set $\cG(N, P;S,X)$ in $\overline{\Teich(S)}$ is discrete.
\label{discrete}
\end{lemma}

{\bf Proof:}
Assume that $Y^{\tau_k}_\geod \rightarrow Y^\tau_\geod$ is a convergent sequence in $\cG(N, P;S,X)$.  Then we can choose an $N > 0$ such that $d_{\rm WP}(Y^{\tau_k}_\geod, Y^{\tau}_\geod)< \delta_0/2$ for $k > N$ where $\delta_0$ is the constant in Wolpert's strata separation  theorem (Theorem \ref{strata}). By the triangle inequality we also have $d_{\rm WP}(Y^{\tau_k}_\geod, Y^{\tau_l}_\geod)< \delta_0$ for $k, l > N$. Thus by Wolpert's strata separation  theorem 
 we have $i(\tau_k, \tau_l) = i(\tau_k,\tau) = 0$ for $k,l > N$. This implies that $\tau_k$ can be only a finite number of possibilities for $k > N$ and therefore $\cG(N, P;S,X)$ is discrete. 
\eproof

We will also be interested in the manifold obtained by {\em drilling}
the curves in $\tau$ from the interior of $N$. We set notation here: 
\begin{itemize}

\item Set $W \cong \partial N \times [0,1]$ be a collar neighborhood of
$\partial N$ with $\partial_0 W  = \bdry N \times \{0\}$ the component of
the boundary lying in $\interior (N)$.

\item Set $\tau_0 = \tau \times \{0\}$ be copies of $\tau$ isotoped
into $\interior(N)$, lying on $\partial_0 W$. 

\item Let $\hat N$ be the compact 3-manifold obtained removing open
tubular neighborhoods $\mathcal{N} (\tau_0)$  of $\tau_0$.

\item Note that $\partial \hat N$ is the union of $\partial N$ and a
torus for each component of $\tau_0$. Let $\hat P$ be the union of $P$
and the new tori 
in $\partial \hat N$ so there is a natural homeomorphism from
$\partial N - P$ to $\partial \hat N - \hat P$. 
\end{itemize}

There is an inclusion  $\iota\colon\hat N \hookrightarrow N$
that
restricts to a homeomorphism from $\partial\hat N - \hat P$ to $\partial N - P$. Therefore $MP(\hat N, \hat P; S, X)$ is also parameterized by $\Teich(S)$.
\begin{itemize}
\item Given $Y \in \Teich(S)$, $\hat M_Y \in MP(\hat N, P;S,X)$ is the hyperbolic manifold such that $\iota$ extends to a conformal map  between the conformal boundary of $\hat M_Y$ and $M_Y$.
\item $\hat\phi_Y$ is the Schwarzian quadratic differential  for the projective structure on $Y$ induced  by $\hat M_Y$.
\end{itemize}

%% When the curves $\tau$ are isotoped inside $N$ to $\hat N$ they are
%% contained in a collar neighborhood of $S$ in $N$. Therefore

There is a natural embedding 
$$j\colon N\to \hat N$$
obtained by including the submanifold $N \setminus \interior
(W \cup \mathcal{N}(\tau_0)) \includesin N$
such that the composition $j\circ \iota$ is isotopic to the identity
and $j$ is a homeomorphism from $\partial N - (P \cup S)$ to
$\partial \hat N - (\hat P \cup S)$.

For every hyperbolic manifold in $MP(\hat N, \hat P; S, X)$ this
embedding induces a cover that lies in $MP(N, P_\tau; S_\tau,
X)$. That is, there is an induced map
$$j^*\colon MP(\hat N, \hat P; S, X) \to MP(N, P_\tau; S_\tau, X)$$
between the deformation spaces and we set
$$M_{\hat Y} = j^*(\hat M_Y).$$

%% file: smallL2.tex
%% !TEX root =renorm_wp.tex

\subsubsection*{Outline of the proof of Theorem \ref{small-L2}}
If $\|\phi_Y\|_\infty$ is small the proof is straightforward:
Thurston's skinning map is a map from $MP(N, P;S, X)$ to itself that
has a fixed point at the totally  geodesic structure. By a theorem of
McMullen the skinning map is contracting and therefore we obtain a
bound on the distance from $Y$ to $Y_\geod$ if we can bound distance
between $Y$ and its first skinning iterate. When $\|\phi_Y\|_\infty$
is small a classical result of Ahlfors-Weill bounds this initial
distance. 

A key element of our investigation involves understanding the behavior
of the $L^\infty$-norm when the $L^2$-norm is small. 
In particular, the pointwise-norm of
$\phi_Y$ may be large in the thin parts of $Y$ which we will need to
pinch to nodes. There are several steps to the proof: 
\begin{itemize}

\item We choose $\tau$ to be the simplex of short curves on $Y$.
 A version of the {\em drilling theorem} bounds the
$L^2$-norm of $\phi_Y -\hat\phi_Y$ in terms of the length of $\tau$. We use this to bound the pointwise norm of $\hat\phi_Y$ outside of the standard collars of $\tau$.

\item 
Using the above bullet and a modification of some classical arguments, this bounds $\|\phi_{\hat Y}\|_\infty$. We are then in position to use McMullen's contraction theorem to bound the distance between $\hat Y$ and $Y^\tau_\geod$.

\item We also have that $Y- \tau$ conformally embeds in $\hat Y$ which implies that $Y$ and $\hat Y$ are close in the Weil-Petersson completion. Together, this and the previous bullet imply the theorem.
\end{itemize}

\subsection{Choosing the curves  to drill}
\newcommand{\cd}{{c_{\rm drill}}}
\newcommand{\ld}{{\ell_{\rm drill}}}

As we noted in the outline, a bound in $\|\phi_Y\|_2$ does not give a bound on $\|\phi_Y\|_\infty$. However, we have the following bound on the pointwise norm that depends on the injectivity radius.  For $Y$ a hyperbolic surface and $z \in Y$ we define  ${\rm inj}(z)$ to be the injectivity radius of $z$ in the hyperbolic metric on $Y$. For simplicity, we define the truncated injectivity radius by ${\rm inj}_Y^-(z) = \min\{{\rm inj}_Y(z), \epsilon_2\}$ where $\epsilon_2 = \sinh^{-1}(1)$ is the Margulis constant in dimension $2$.
\begin{prop}[Bridgeman-Wu, \cite{bridgeman2019uniform}]\label{pointbound}
Let $\phi \in Q(Y)$ then 
$$\|\phi(z)\| \le \frac{\|\phi\|_2}{\sqrt{{\rm inj}_Y^-(z)}}.$$
\end{prop}

As a first step we show that after an appropriate choice for $\tau$, we can obtain a pointwise bound on $\hat\phi_Y$ outside of the standard collars of $\tau$. For this we will need the following bound on the $L^2$-norm.
\begin{theorem}[Bridgeman-Bromberg, \cite{BBcone}]\label{conedrill}
There exists constants  $\cd, \ld > 0$ with $\ld <2\epsilon_2$ such that the following holds. Given $Y \in \Teich(S)$ and a simplex  $\tau$ in $\mathcal C(S)$ such that for all $\beta\in \tau$, $\ell_\beta(Y)\le \ld$ then
$$\|\phi_Y-\hat\phi_Y\|_2 \le \cd\sqrt{\ell_\tau(Y)}$$
where
$\ell_\beta(Y)$ is the length of $\beta$ in $Y$.
\end{theorem}

\subsubsection*{Fixing a universal constant}

We first prove that we can choose the simplex $\tau$ such that $||\hat\phi_Y(z)||$ is small for $z \in Y$ in the complement of the standard collars of $\tau$.

\begin{theorem}{}\label{drillL2bound}
 Assume that $Y\in\Teich(S)$ with $\|\phi_Y\|_2^{\frac2{2n(S)+3}} \le \ld$. There exists an $\ell= \ell(Y)>0$ with
$$\ell \le \|\phi_Y\|_2^{\frac2{2n(S)+3}}$$
such that the following holds. Let $\tau$ be the simplex in $\mathcal C(S)$ of all curves with length $\le\ell$. Then for $z\in Y$  in the complement of the standard collars of $\tau$
$$\|\hat\phi_Y(z)\| \le C_0\sqrt{n(S)} \|\phi_Y\|_2^{\frac{2}{2n(S)+3}}.$$
for $C_0 = \sqrt{2}(\cd+1)$.\end{theorem}

{\bf Proof:} 
Let $\Lambda = \|\phi_Y\|_2^{\frac{2}{2n(S)+3}}$ and let $\ell_k = \Lambda^{2k+1}$. As $\Lambda <2\epsilon_2$ there are at most $n(S)$ curves of length $\le \Lambda$ so there must be some integer $k$ with $0 \le k \le n(S)$ such that $Y$ has no curves of length in the interval $(\ell_{k+1}, \ell_{k}]$. Let $\ell = \ell_{k+1} \le \ell_0 = \|\phi_Y\|_2^{\frac2{2n(S)+3}}$ and let $\tau$ be the simplex in $\cC(S)$ of all curves of length $\le \ell$ on $Y$. 

By Theorem \ref{conedrill} we have
$$\|\phi_Y - \hat\phi_Y\|_2 \le \cd\sqrt{\ell_\tau(Y)}.$$
As $\ell_\tau(Y) \le n(S) \Lambda^{2k+3}$ we have
$$\|\hat\phi_Y\|_2 \leq \|\phi_Y\|_2 + \|\phi_Y - \hat\phi_Y\|_2 \le \Lambda^{n(S)+\frac32} + \cd\sqrt{n(S)}\Lambda^{k+\frac32}.$$
As $Y$ contains no curves of length in the interval $(\ell_{k+1}, \ell_{k}]$ every point in the complement of the standard collars of $\tau$ has injectivity radius $> \ell_{k}/2 = \Lambda^{2k+1}/2$. Therefore if $z\in Y$ is in the complement of the standard collars of $\mathcal C$ then by Proposition \ref{pointbound}
\begin{eqnarray*}
\|\hat\phi_Y(z)\| &\le& \frac{\|\hat\phi_Y\|_2}{\sqrt{\ell_{k}/2}}\\
& \le & \frac{\Lambda^{n(S)+\frac32} +\cd\sqrt{n(S)}\Lambda^{k+\frac32}}{\sqrt{\Lambda^{2k+1}/2}} \\
& \le & \sqrt2\left(\Lambda + \cd\sqrt{n(S)} \Lambda\right)\\
&\le& \sqrt2\left(1 + \cd\sqrt{n(S)}\right) \Lambda\\
&\le& C_0\sqrt{n(S)} \Lambda
\end{eqnarray*}
where $C_0 = \sqrt{2}(1+\cd)$ is a universal constant. \eproof

We can now prove:
\begin{theorem}\label{progress}
If $Y\in\Teich(S)$ with $\|\phi_Y\|_2^{\frac2{2n(S)+3}} \le \min\{\ld,2\sinh^{-1}(1/2)\}$ then there is a simplex $\tau \in\cC(S)$ and a  $\hat Y \in \Teich(S_\tau) \subseteq \overline{\Teich(S)}$ such that 
\begin{enumerate}[(a)]
\item $d_{\rm WP}(Y,\hat Y) \le \frac{2\pi}{\sqrt{\sinh^{-1}(1/2)}} \sqrt{n(S)} \|\phi_Y\|_2^{\frac1{2n(S)+3}}$

\item $\|\phi_{\hat Y}\|_\infty \le C_1\sqrt{n(S)}\|\phi_Y\|_2^{\frac1{2n(S)+3}}$
\end{enumerate}
where $C_1 = 9\sqrt{2}(C_0+1)$. 

\end{theorem}
{\bf Proof of Theorem \ref{progress} (a):}  Let $\tau$ be the simplex given by Theorem \ref{drillL2bound} and $\hat Y \in \Teich(S_\tau)$ the surface with $j^*(M_{\hat Y}) = \hat M_Y$.
We need a few topological observations about the cover.
Let $\hat N_Y$ be the union of $\hat M_Y$ and the conformal boundary $\partial_c \hat M_Y$. Similarly define $N_{\hat Y}$ for $M_{\hat Y}$. While $M_{\hat Y}$ covers $\hat M_Y$, the associated cover $\check N_{\hat Y}$ of $\hat N_Y$ is a proper submanifold of $N_{\hat Y}$. In particular the boundary of $\check N_{\hat Y}$ will be a union of covers of components of $\partial_c \hat M_Y$. The restriction of $j$ to $S\backslash \tau$ will be homotopic to an embedding of $S\backslash \tau$ into $Y \subset \partial_c \hat M_Y$ and therefore the, possibly disconnected, cover $\check Z$ of $Y$ associated to $S\backslash \tau$ {will be components of} $\partial \check N_{\hat Y}$. In particular, $\check Z$ will conformally embed in $\hat Y$.
By assumption, $\|\phi_Y\|_2^{\frac2{2n(S)+3}} \le 2\sinh^{-1}(1/2)$ and therefore, by Theorem \ref{drillL2bound}, for each $\beta\in\tau$,
$$\ell_\beta(Y) \le 2\sinh^{-1}(1/2)$$
so we can apply Theorem \ref{wp estimate} to get
$$d_{\rm WP}(Y, \hat Y) \le \frac{2\pi}{\sqrt{\sinh^{-1}(1/2)}}\sqrt{\ell_\tau(Y)} \leq \frac{2\pi}{\sqrt{\sinh^{-1}(1/2)}}\sqrt{ n(S)} \|\phi_Y\|_2^\frac{1}{2n(S)+3}.$$
 \eproof

To obtain our bound on $||\phi_{\hat Y}||_\infty$ we will first need the following generalization of  the Kraus-Nehari bound on the norm of the Schwarzian.
\begin{lemma}
\label{bigdisk}
Let $f\colon\Delta\to\Delta$ be univalent and assume that for $z \in \Delta$ the image $f(\Delta)$ contains a hyperbolic disk of radius $r$ centered at $z$. Then $\|Sf(z)\| \le \frac32\sech \frac{r}{2}$.
\end{lemma}

{\bf Proof:} The proof  is a refinement of the classical proof of the Kraus-Nehari theorem.

Assume that $z = f(z)= 0$. By applying the Schwarz Lemma to the restriction of $f^{-1}$ to the hyperbolic disk of radius $r$ we see that $|f'(0)| \ge \tanh \frac{r}{2}$. If we let $g(z)= f'(0)/f(1/z)$ we have the expansion
$$g(z) = z +\sum_{n=0}^\infty b_nz^{-n}.$$
Note that the domain of $g$ is $\{z\in \hat{\mathbb C}\ | \ |z|>1\}$ and that $|g(z)| > \tanh \frac{r}{2}$ for $z$ in the domain. As in the proof of Nehari's theorem we can also calculate to see that $Sf(0)=-6b_1$. As  the conformal factor for the area form of the hyperbolic metric on $\Delta$ at $z=0$ is $4$, we obtain $\|Sf(0)\| =\frac32|b_1|$.  Let $C_\rho$ be the circle of radius $\rho$ centered at $0$ with $\rho>1$. Then the Euclidean area $m_\rho$ in $\mathbb C$ bounded by $g(C_\rho)$ is given by
$$m_\rho = \pi\rho^2-\pi\sum_{n=1}^\infty n|b_n|^2\rho^{-2n}.$$
Since, for all $\rho>1$, $C_\rho$ will contain the disk of radius $\tanh \frac{r}{2}$ centered at $0$ we have that $m_\rho > \pi \tanh^2 \frac{r}{2}$ and by letting  $\rho\to 1$ we have
$$\pi \tanh^2 \frac{r}{2} \le \pi - \pi \sum_{n=1}^\infty n|b_n|^2 \le \pi - \pi |b_1|^2.$$
The estimate follows.
\eproof

{\bf Proof of Theorem \ref{progress} (b):} Choose $\epsilon$ such that
$$ \|\phi_Y\|_2^{\frac2{2n(S)+3}} = 2\epsilon$$
and note that $2\epsilon \le \epsilon_2$.

To simplify notation will assume that $S$ contains a single component. The general case easily follows.
The hyperbolic 3-manifolds $\hat M_Y$ and $M_{\hat Y}$ are uniformized by Kleinian groups which we denote $\hat \Gamma_Y$ and $\Gamma_{\hat Y}$ where $\Gamma_{\hat Y}$ is a subgroup of $\hat \Gamma_Y$. Fix a component $\widehat W$ of $\hat Y$. Let $\Omega_Y$ be a component of the domain of discontinuity of $\hat \Gamma_Y$ that covers $Y$ in $\partial \hat N_Y$ and $\Omega_{\widehat{W}}$  a component of the domain of discontinuity of $\Gamma_{\hat{Y}}$ that covers $\widehat{W}$. We can further assume that $\Omega_Y \subset \Omega_{\widehat{W}}$ and that the stabilizer $\Gamma_{\widehat{W}}$ of $\Omega_{\widehat{W}}$ in $\Gamma_{\hat Y}$ also stabilizes $\Omega_Y$. In particular the quotient $\Omega_{\widehat{W}}/\Gamma_{\widehat{W}}$ is $\widehat{W}$ while the quotient $\check W = \Omega_Y/\Gamma_{\widehat{W}}$ is a component of the cover $\check Z$ of $Y$ discussed in part (a).

Let $W$ be the convex core of $\check W$. Then $W$ is bounded by a collection of geodesics that map to components of the geodesic representatives of $\tau$ in $Y$ under the covering map $\check W\to Y$. Also observe that the interior of $W$ embeds in $Y$ under this covering. As the length of these geodesics is $\le 2\epsilon$ the injectivity radius of $\partial W$ is $\le \epsilon$. Let $\widehat W^{\epsilon_2}$ and $\widehat W^\epsilon$ be the complement of the $\epsilon_2$ and $\epsilon$-cuspidal thin parts of $\widehat{W}$, respectively. 
By the Schwarz Lemma the embedding $\check W\hookrightarrow \widehat{W}$ is a contraction from the complete hyperbolic metric on $\check W$ (which is lifted from $Y$) to the complete hyperbolic metric on $\widehat{W}$. Therefore under this embedding the boundary of $W$ will map into the $\epsilon$-cuspidal thin part of $\widehat{W}$ so $W \supset \widehat W^\epsilon$. As $\widehat W^{\epsilon_2} \subset  \widehat W^\epsilon$ we also have $\widehat W^{\epsilon_2} \subset W$ and therefore $\widehat W^{\epsilon_2}$ embeds in $Y$. Furthermore the standard collars of $\partial W$ in $\check W$ have injectivity radius $\le \epsilon_2$ so they lie in the $\epsilon_2$-cuspidal thin part of $\widehat{W}$ and therefore $\widehat W^{\epsilon_2}$ is in the complement of the standard collars of $\tau$ in $Y$.

Now let $f_Y \colon \Delta\to \Omega_Y$ and $f_{\widehat{W}}\colon\Delta\to\Omega_{\widehat{W}}$ be the associated uniformizing maps.
As  $\Omega_Y \subset \Omega_{\widehat{W}}$ we can factor $f_Y$ through a map $g\colon\Delta\to \Delta$ such that $f_Y = f_{\widehat{W}}\circ g$. 

By (\cite[Lemma 4.5]{Barrett:Diller_contract}) the norm of a quadratic differential achieves it maximum in the complement of the standard neighborhood of the cusps. Therefore to bound $\|\phi_{\hat Y}\|_\infty$ it suffices to bound  $\|\phi_{\widehat{W}}(z)\|$ for $z\in \widehat W^{\epsilon_2}$. 

After fixing a $z \in \widehat W^{\epsilon_2}$ it will be convenient to normalize our uniformizing maps such that $g(0) = 0$ and that $0$ maps to $z$ under the quotient maps to $\widehat{W}$ and $Y$. Then
$$\|\phi_{\hat Y}(z)\| = 4|Sf_{\widehat{W}}(0)| \mbox{ and }  \|\hat\phi_Y(z)\|= 4|Sf_Y(0)|.$$

By the composition rule for Schwarzian derivatives we have
$$Sf_Y(0) = Sf_{\widehat{W}}(g(0)) g'(0)^2 + Sg(0)$$
and therefore (assuming that $g(0) = 0$)
$$\|\phi_{\hat Y}(z)\|=\|Sf_{\widehat{W}}(0)\| \le (\|Sf_Y(0)\| + \|Sg(0)\|)/|g'(0)|^2.$$
We now need to bound the individual terms on the right. 

As $ \widehat W^{\epsilon_2}$ is in the complement of the standard collars of $\tau$ in $Y$, by Theorem \ref{drillL2bound}
$$\|Sf_Y(0)\| = \|\phi_Y(0)\| \le 2C_0\sqrt{n(S)}\epsilon.$$

We would like to apply Lemma \ref{bigdisk} to bound $\|Sg(0)\|$ but to do so we need to bound from below the distance from $0$ to $\Delta\backslash g(\Delta)$ in the hyperbolic metric on $\Delta$. This distance is bounded below by the distance from $\widehat W^{\epsilon_2}$ to $\widehat{W}\backslash \check W$ in the hyperbolic metric on $\widehat{W}$ and this distance in turn is bounded below by the distance from $\widehat W^{\epsilon_2}$ to $\widehat{W}\backslash W$. To bound this we use that $\widehat W^\epsilon\subset W$ and a simple calculation shows that if $r$ is the distance from $\partial \widehat W^{\epsilon_2}$ to $\partial \widehat W^\epsilon$  then
$$e^r > \frac{\sinh(\epsilon_2/2)}{\sinh(\epsilon/2)} > \frac{\epsilon_2}{\epsilon} \geq 2.$$
The hyperbolic disk of radius $r$ centered at $0$ will be contained in $g(\Delta)$
and Lemma \ref{bigdisk} plus the above bound implies that
$$\|Sg(z)\| \le \frac32 \sech \frac{r}{2} < 3e^{-\frac{r}{2}} < 3\sqrt{\frac{\epsilon}{\epsilon_2}}.$$

Finally we need to bound from below $|g'(0)|$. As in the proof of Lemma \ref{bigdisk} we have $|g'(0)| \ge \tanh \frac{r}{2}$ and given our above bound on $r$ this becomes
$$|g'(0)| \ge \tanh \frac{r}{2} \ge  \frac{1 - \frac{\epsilon}{\epsilon_2}}{1 + \frac{\epsilon}{\epsilon_2}} \geq \frac{1}{3}.$$

Combining our estimates we have
\begin{eqnarray*}
\|Sf_{\widehat{W}}(0)\| \le 9\left(2C_0\sqrt{n(S)} \epsilon +3\sqrt{\frac{\epsilon}{\epsilon_2}}\right) \leq  9\sqrt{2}(C_0+1)\sqrt{n(S)} \|\phi_Y\|_2^{\frac1{2n(S)+3}}.
\end{eqnarray*}
Therefore we let $C_1 = 9\sqrt{2}(C_0+1)$ and the result follows.
\eproof

\newcommand{\ctm}{{A}}
\subsection{Bounds on iteration of the skinning map}
Let $(N,P;S)$ be a relatively acylindrical triple. Then for $Y \in \Teich(S) \cong MP(N,P;S,X)$ we need to show that if $\|\phi_Y\|_\infty$ is small then $d_{\rm WP}(Y, Y_{\rm geod})$ is small.
When $(N,P)$ is acylindrical the proof is a straightforward application of a classical bound of Ahlfors-Weill plus McMullen's contraction theorem for the skinning map. However, in the relatively acylindrical case we will need a slight extension of McMullen's original statement.

The {\em skinning map} 
$$\sigma \colon MP(N,P;S,X) \simeq \Teich(S)  \to \Teich(S)$$
is defined as follows: for each $Y \in \T(S)$, the cover of $M_Y \in MP(N,P;S,X)$ associated to the subgroup $\pi_1(Y) \subset \pi_1(M_Y)$
under inclusion will be quasifuchsian. (If $Y$ is disconnected then
the cover will also be a finite collection of a quasifuchsian
manifolds.) For each connected component of $\partial M_Y$, one component of
the conformal boundary restricts to a homeomorphism to $Y$ under the
covering projection. The other component will be $\sigma(Y)$, the
image of the skinning map for that component. Note that $Z\in \Teich(S)$ is in $\mathcal G(N,P;S,X)$ if and only if $Z$ is a fixed point for $\sigma$.

The skinning map is a smooth map and we will be interested in bounding its derivative so that we can apply the contraction mapping principle. The estimate we need from McMullen essentially works as written in \cite{McMullen:iter} but there a few differences in the relative case that we highlight. Given $Y\in \Teich(S)$ let $\Gamma$ be the Kleinian group that uniformizes $M_Y \in MP(N,P;S,X)$ and let $\Omega$ be the domain of discontinuity of $\Gamma$. If the pair $(N,P)$ was acylindrical then every component of $\Omega$ would be a Jordan domain and the stabilizer of every component would be a quasifuchsian group. Furthermore if $D_0$ and $D_1$ are distinct components of $\Omega$ then either their closures are disjoint and the intersection of their stabilizers is trivial or the intersection is a point and the intersection of their stabilizers is an infinite cyclic group generated by a parabolic. In the  relatively acylindrical this will not hold. However, if we let $\Omega_Y$ be those components of $\Omega$ that cover $Y$ then these properties do hold for the components in $\Omega_Y$. The second key point is that a tangent vector of $MP(N, P; S,X)$ is represented by a $\Gamma$-invariant Beltrami differential $\mu$ that is supported on $\Omega_Y$. With these two observations one sees that McMullen's proof in the acylindrical case extends to the relatively acylindrical case:
\begin{theorem}[McMullen, {\cite[Thm 6.1, Cor 6.2]{McMullen:iter}}]
If $(N,P;S,X)$ is relatively acylindrical then for $Y \in \Teich(S)$ 
$$||d\sigma_Y||_\infty \le \lambda(S)<1$$
where $\lambda(S)$ depends only on the topology of $S$.
\label{acyl_contraction}
\end{theorem}

The contraction mapping principle implies that $\sigma^n(Y) \to Z$ with $\sigma(Z) = Z$ and
$$d_{\Teich}(Y,Z) \le \frac{d_{\Teich}(Y, \sigma(Y))}{1-\lambda(S)}.$$
To complete the proof of Theorem \ref{geodesic distance} we need to bound $d(Y, \sigma(Y))$. This is a direct consequence of the Ahlfors-Weill quasiconformal reflection theorem:
\begin{theorem}[Ahlfors-Weill, {\cite[Theorem 5.1]{Lehto:book:univalent}}]
Let $Y \in  \Teich(S)$ and $\phi_Y$ be the associated quadratic differential on $Y$. If $||\phi_Y||_\infty < 1/2$ then
$$d_{\Teich}(Y,\sigma(Y)) \leq \frac12 \log\frac{1+2 ||\phi_Y||_\infty}{1-2\|\phi_Y\|_\infty}.$$
\label{AhlforsWeill}
\end{theorem}
If $\|\phi_Y\|_\infty \le \frac13$ then an easy estimate of the right hand side gives that 
$$d_{\Teich}(Y,\sigma(Y)) \le  3 \|\phi_Y\|_\infty$$
and therefore
$$d_{\Teich}(Y,Z) \le \frac{3}{1-\lambda(S)}\|\phi_Y\|_\infty.$$
By a result of Linch (see \cite{Linch:wp}), $d_{\rm WP} \le \sqrt{\area(S)} d_{\Teich}$ and we have the following result:
\begin{theorem}\label{geodesic distance}
Let $(N,P;S)$ be relatively acylindrical.
Then for all $Y\in \Teich(S)$ with $\|\phi_Y\|_\infty \le 1/3$ we have
$$\frac{d_{\rm WP}(Y,Z)}{\sqrt{\area(Y)}} \le d_{\Teich}(Y,Y_{\rm geod}) \le \frac{3\|\phi_Y\|_\infty}{1-\lambda(S)}$$
where $\lambda(S)$ is the contraction constant from Theorem \ref{acyl_contraction}.
\end{theorem}

{\bf Remark.} McMullen's proof is not effective and this is the one place in our proof we were don't control the growth rate of the constants in terms of genus. However, we have made some effort to isolate this from the constants that we do control. 

\subsection{Proof of Theorem \ref{small-L2}}
We now put together the results above. We first restate Theorem \ref{small-L2}, but here we carefully control the constants.

\begin{theorem}\label{nearnode}
There are a universal constants $K_0$ and $\epsilon_0$ such that if 
$$A(\epsilon, S)  = \left(\frac{K_0\epsilon (1-\lambda(S))}{n(S)}\right)^{2n(S)+3}$$
and $Y \in \Teich(S)$ with $\|\phi_Y\|_2 \leq A(\epsilon, S)$ and $\epsilon\leq \epsilon_0$ then there exists $Y^\tau_\geod \in \mathcal G$ with $d_{\rm WP}(Y, Y^\tau_\geod) \leq \epsilon$.
\end{theorem}

{\bf Proof:} By Theorem \ref{progress}, there are universal constants $\ld, C_1 > 0$ such that if $||\phi_Y||_2^{\frac{2}{2n(S) + 3}} \leq \ld$ then there is a simplex $\tau$ in $\mathcal C(S)$ such that after drilling curves $\mathcal C$, 
$$||\phi_{\hat Y}||_\infty \leq C_1\sqrt{n(S)} ||\phi_Y||_2^{\frac{1}{2n(S)+3}} \qquad d_{\rm WP}(Y,\hat Y) \leq \frac{2\pi}{\sqrt{\sinh^{-1}(1/2)}}\sqrt{n(S)} ||\phi_Y||_2^{\frac{1}{2n(S)+3}}.$$
Assuming that $\|\phi_{\hat Y}\|_\infty \le 1/3$ we can apply Theorem \ref{geodesic distance} to $(N_\tau, P_\tau; S_\tau)$ to see that 
$$d_{\rm WP}(\hat Y, Y^\tau_\geod) \leq \frac{3\sqrt{\area(\hat Y)}}{1-\lambda(S)} \|\phi_{\hat Y}\|_\infty \le \frac{2\sqrt{3\pi}C_1n(S)}{1-\lambda(S)} ||\phi_Y||_2^{\frac{1}{2n(S)+3}}$$
since $\area(\hat Y) = \area(Y) = 4\pi n(S)/3$.
Then by the triangle inequality and the fact that $C_1 > 1$ we have
$$d_{\rm WP}(Y, Y^\tau_\geod) \leq d_{\rm WP}(Y,\hat Y)+d_{\rm WP}(\hat Y,Y^\tau_\geod)\le \frac{4\sqrt{3\pi}C_1n(S)}{1-\lambda(S)} ||\phi_Y||_2^{\frac{1}{2n(S)+3}}$$
We let $K_0= 1/(4\sqrt{3\pi}C_1).$  Recounting our progress, if
$$\|\phi_Y\|_2 \le A(\epsilon,S) =   \left(\frac{K_0\epsilon(1-\lambda(S))} {n(S)}\right)^{2n(S)+3}$$
we have
$$d_{\rm WP}(Y, Y^\tau_\geod) \le \epsilon$$
assuming that $\|\phi_Y\|_2^{\frac{2}{2n(S)+3}}  < \ld$ and $\|\phi_{\hat Y}\|_\infty \le 1/3$. However, if we let 
$$\epsilon_0 = \min\left\{\frac{\sqrt{\ld}}{K_0},  4\sqrt{\frac{\pi}{3}}\right\}$$
and $\epsilon < \epsilon_0$ then
$$\|\phi_Y\|_2^{\frac{2}{2n(S)+3}} \le  \left(\frac{K_0\epsilon(1-\lambda(S))} {n(S)}\right)^2 \le (K_0 \epsilon)^2 < \ld$$
and
\begin{eqnarray*}
\|\phi_{\hat Y}\|_\infty &\le& C_1\sqrt{n(S)} ||\phi_Y||_2^{\frac{1}{2n(S)+3}}\\
& \le & C_1\sqrt{n(S)}\frac{K_0\epsilon(1-\lambda(S))} {n(S)}\\
& \le & C_1K_0\epsilon = \frac{\epsilon}{4\sqrt{3\pi}}  \le \frac1 3.
\end{eqnarray*}
This completes the proof.
 \eproof

\subsubsection*{Hyperbolic manifolds with cylinders and compression disks}
We conclude this section with a discussion of where we use the relative acylindricity of $(N,P;S)$. For simplicity, in this discussion we will assume that both $P$ and $S$ are empty.

The first problem that can occur is in Theorem \ref{drillL2bound} and its application. In particular, it can happen that two or more curves in $\tau$ may be homotopic in $N$ or even homotopically trivial in $N$ if $N$ has compressible boundary. In this case, the manifold $\hat N$ will not be hyperbolizable. If $N$ has incompressible boundary, this problem can be corrected by only removing a single curve from $N$ for each homotopy class (in $N$) of curves in $\tau$. With this change, Theorem \ref{drillL2bound} we still hold but we cannot define the embedding of $N$ in $\hat N$ and therefore cannot carry through the proof of Theorem \ref{progress}.

If none of the curves in $\tau$ are homotopic in $N$ then the proof up to and including Theorem \ref{progress} go through. However, if the pared manifold $(N, P_\tau)$ is not acylindrical then Theorem \ref{acyl_contraction}, McMullen's contraction theorem, will fail. In fact, the deformation space $MP(N, P_\tau)$ contains a hyperbolic structure whose convex core boundary is totally geodesic if and only if $(N,P_\tau)$ is acylindrical or is a pared $I$-bundle. 

We expect that the only problem that can occur is the first one. We have the following conjecture.
\begin{conj}\label{general_drill}
Let $M$ be a convex co-compact hyperbolic 3-manifold with $\phi$ the Schwarzian quadratic differential for the projective boundary of $M$. If $\|\phi\|_2$ is small then either:
\begin{itemize}
\item There exists a geometrically finite structure $M'$ on $N$ with totally geodesic convex core boundary and $d_{\rm WP}(\partial_c M, \partial_c M')$ is small.

\item There are two or more short curves on $\partial_c M$ that are homotopic in $M$.
\end{itemize}
\end{conj}
In particular, if no two curves in $\tau$ are homotopic in $M$ we expect that $(N,P_\tau)$ is an acylindrical pair even when $N$ itself is not acylindrical.

%% file: bersslice.tex
%% !TEX root =renorm_wp.tex

\section{W-volume and renormalized volume}
Given a convex submanifold $N$ with smooth boundary such that $N \hookrightarrow M$ is a homotopy equivalence, the {\em W-volume} of $N$ is defined to be
$$W(N) = \vol(N) -\frac{1}{2}\int_{\partial N} HdA$$
where $H$ is the mean curvature of $\partial N$. 

The $W$ volume has many nice analytic properties that make it a useful tool for studying hyperbolic manifolds. We let $N_t$ be the $t$-neighborhood of $N$. The nearest point retraction from $M$ to each $N_t$ extends to a diffeomorphism from $\partial_c M$ to $\partial N_t$ and using this retract we pull back the induced metrics on $\partial N_t$ to metrics $I_t$ on $\partial_c M$. Then
$$I^*(x,y) = \lim_{t\rightarrow \infty} \frac{1}{\cosh^2(t)} I_t(x,y)$$
is a well defined metric in the conformal class of $\partial_c M$.

For $N \subset M$ we will denote by $\rho_N$ the metric at infinity on $\partial_cM$. The $W$ volume has the following properties.

\begin{prop}{(Krasnov-Schenker, \cite{KS08})}\label{scaling properties}
Let $N\subset M$ be a compact, convex submanifold of a convex cocompact hyperbolic 3-manifold $M$ and let $N_t$ be the $t$-neighborhood of $N$. Then
\begin{enumerate}
\item The metric $\rho_N$ is in the conformal class of $\partial_c M$;

\item $\rho_{N_t}  = e^{2t}\rho_N$;

\item $W(N_t) = W(N) -t\pi\chi(\partial N)$.
\end{enumerate}
Furthermore if $\rho$ is any smooth conformal metric on $\partial_c M$ then for $t$ sufficiently large there exists a convex submanifold $X_t \subset M$ with $\rho_{X_t} = e^{2t}\rho$.
\end{prop}
Using this proposition,  the $W$-volume of any smooth conformal metric $\rho$ on $\partial_c M$ is defined by
$$W(\rho) = W(N_t(\rho)) + t\pi\chi(\partial M)$$
for $t$ sufficiently large. The proposition above implies that $W(\rho)$ doesn't depend on the choice of $t$. With this setup we can now define the {\em renormalized volume} $V_R$ by setting
$$V_R(M) = W(\rho_M)$$
where $\rho_M$ is the unique hyperbolic metric on $\partial_c M$.

\subsubsection*{Convex cores}
Perhaps the most natural convex submanifold of a convex co-compact hyperbolic 3-manifold $M$ is the {\em convex core} $C(M)$. The boundary of the convex core is not in general smooth, so we cannot use the previous definition to define the $W$-volume of $C(M)$. However, there is a natural way to extend $W$-volume to this setting (see the discussion in \cite{BBB}) and for the convex core we have
$$W(C(M)) = V_C(M) - \frac{1}{4}L(\beta_M)$$
where $\beta_M$ is the {\em bending lamination} of the boundary of the convex core and $L(\beta_M)$ is its length (as a measured lamination). The convex core also induces  a natural metric at infinity, called the {\em projective metric} (so called as Thurston gave a definition that is intrinsic to the induced projective structure on $\partial_c M$). We will be interested in a hybrid metric that is the hyperbolic metric on some components of $\partial_c M$ and the projective metric on the others. We have the following:
\begin{prop}\label{Wbound}
Let $M$ be a convex co-compact hyperbolic 3-manifold and $\partial_c M = X\sqcup Y$,  a disjoint union of connected components of $\partial_c M$. Let $\sigma$ be the hyperbolic metric on $X$ and the projective metric on $Y$. Let $\beta_Y$ be the bending lamination of the components of the boundary of $C(M)$ that faces $Y$. Then
$$W(\sigma) - \frac14 L(\beta_{Y}) \le V_R(M) \le W(\sigma).$$
In particular if $Y=\partial N$ we have
$$V_C(M) - \frac12 L(\beta_M) \le V_R(M) \le V_C(M) -\frac14 L(\beta_M).$$
\end{prop}
By the definition of the $W$-volume of the convex core, the two statements are equivalent for the case $X=\emptyset$ and was proven in \cite[Theorem 3.7]{BBB}. Furthermore,  the proof trivially extends to the relative case above.

\section{The variational formula}
Recall that if $(N;S)$ is a pair such that each component of $S$ is incompressible in $N$ then $MP(N:S,X)$ is parameterized by $\Teich(S)$ and therefore we can view renormalized volume as a function
$$V_R\colon \Teich(S) \to \mathbb R.$$
We recall the variational formula:

\medskip
\noindent
{\bf Theorem \ref{variational}}
{\em Given $Y \in \Teich(S)$ and $\mu \in T_Y\Teich(S)$ we have
$$dV_R(\mu) = \Re \int_{\partial_c M_Y} \phi_Y \mu.$$
}

Therefore the Weil-Petersson gradient of $V_R$ has norm $\|\phi_Y\|_2$. By the classical bound of Kraus-Nehari for the Schwarzian of univalent functions we have that $\|\phi_Y\|_\infty \le \frac32$. 
As a corollary we have:
\begin{cor}\label{extension}
The Weil-Petersson norm of the gradient of $V_R$ is bounded by $\frac32 \sqrt{\area(Y)} = \sqrt{3\pi n(S)}$. In particular $V_R$ is Lipschitz with respect to the Weil-Petersson metric and therefore extends to a continuous function on the Weil-Petersson completion.
\end{cor}
Note if $S$ is not incompressible in $N$ then we cannot apply the Kraus-Nehari theorem to bound the norm of the gradient and in fact there is no upper bound of the gradient in this setting.

We now assume that $(N;S)$ is relatively acylindrical and
recall that $\mathcal G = \mathcal G(N;S,X)$ is the collection of $Y\in \overline{\Teich(S)}$ such that the component of the boundary convex core of $M_Y$ facing $Y$ is totally geodesic. 
\begin{prop}\label{no_local_min}
Given $\tau$  a non-empty simplex in $\mathcal C(S)$, let $Y^\tau_\geod$ be the unique surface in $\mathcal G \cap \Teich(S_\tau)$.  Then for $t>0$ there is a one parameter family $Y_t \in \Teich(S)$ with $Y_t \to Y^\tau_\geod$ as $t\to 0$ with $V_R(Y_t) < V_R(Y^\tau_\geod)$.
\end{prop}

{\bf Proof:} 
By a construction of Bonahon-Otal (\cite{Bonahon:Otal:shortgs}) there exists a one-parameter family $M_\theta \in MP(N;S,X)$ where the bending lamination $\beta_\theta$ of the components of the convex core facing $S$ have support $\tau$ and bending angle $\theta$. In the parameterization $MP(N;S,X) \cong \Teich(S)$ the manifolds $M_\theta$ correspond to $Z_\theta \in \Teich(S)$. We also let $\sigma_\theta$ be the hybrid metric that is the projective metric on $Z_\theta$ and the hyperbolic metric on $X$. Let $\phi_\theta$ be the Schwarzian quadratic differential on $Z_\theta$.

As part of the construction, Bonahon-Otal show that $M_{Z_\theta}$ converges to $M_{Y^\tau_\geod}$ in the {\em algebraic topology} on $GF(N,S;X)$. Unfortunately what we need is
that $Z_\theta \to Y^\tau_\geod$ in $\overline{\Teich(S)}\cong GF(N,S;X)$ where the topology is the metric topology of the Weil-Petersson completion. These two topologies are not homeomorphic. While the convergence we need could be proven using the notion of strong convergence of Kleinian groups and techniques well-known to experts, we will instead give a proof more in  line with the methods from this paper.

We first note that from the construction it follows that $L(\beta_\theta) \to 0$ as $\theta\to \pi$. In \cite{BBwvol} it is shown that
$$\|\phi_\theta\|_2 \le \frac52 \sqrt{L(\beta_\theta)}$$
and therefore we also have $\|\phi_\theta\|_2 \to 0$ as $\theta\to \pi$.  Theorem \ref{nearnode} then implies that $Z_\theta$ accumulates on $\mathcal G$. As $\mathcal G$ is discrete (see Lemma \ref{discrete}),   $Z_\theta$ must limit to a unique point. It also follows from construction that the length of a curve $\gamma$ on $Z_\theta$ limits to zero if and only if $\gamma$ is in $\tau$ so any limit for $Z_\theta$ will be in the strata for $\tau$. Together this implies that $Z_\theta\to Y^\tau_\geod$.

By Corollary \ref{extension}, $V_R$ extends to a continuous function on $\overline{\Teich(S)}$. Combining this with Proposition \ref{Wbound} and the fact the $L(\beta_\theta) \to 0$ we have
$$\underset{\theta\to\pi}{\lim} W(\sigma_\theta) = \underset{\theta\to\pi}{\lim} V_R(M_\theta) = V_R(Y^\tau_\geod).$$
We will show that
$$V_R(M_\theta) \le W(\sigma_\theta) < V_R(Y^\tau_\geod)$$
which will give the result.

For this we use the variational formula
$$\frac{d}{d\theta} W(\sigma_\theta) = \frac14(\ell(\theta) - \theta \ell'(\theta))$$
where $\ell(\theta)$ is the sum of the length of the curves in $\tau$ on in $M_\theta$. If $X =\emptyset$ then by the Schl\"afli formula
$$\frac{d}{d\theta} V_C(M_\theta) = \frac12 \ell(\theta)$$
and the variational formula follows from differentiating the formula for $W$-volume of the convex core and the noting that $L(\beta_\theta) = \theta\ell(\theta)$. In general, if $\rho_t$ is a family of metrics on $\partial N$ then the variation of $W$-volume will have a term for each component of the boundary and if $\tilde\rho_t$ is a another family of metrics that agrees with $\rho_t$ on a component $S$ of $\partial N$ then the term for both variations on $S$ will be the same. In our case $\sigma_\theta$ is the hyperbolic metric on $X$ for all $\theta$ so the variation of $W$-volume on $X$ is zero. On $Z_\theta$, $\sigma_\theta$ is the projective metric so on $Y$ the variation is the same as the variation of the $W$-volume of the convex core. This gives the variational formula.

We can now complete the proof. By Choi-Series (\cite{choi2006}) $\ell'(\theta) < 0$ which implies that $W(\sigma_\theta) < V_R(Y^\tau_\geod)$. We can also see this directly by integrating to get
$$V_R(Y^\tau_\geod) - W(\sigma_T)  = \frac14 \int_T^\pi \ell(\theta) d\theta + \frac18 T\ell(T)>0.$$
We then define $Y_t$ by reparameterizing $Z_\theta$ via an orientation reversing homeomorphism from $(0,\infty)$ to $(0,\pi)$. Thus we have $V_R(Y_t) < V_R(Y^\tau_\geod)$ as required.
\eproof

%% file: endgame.tex
%% !TEX root =renorm_wp.tex

\section{Lower bounds on renormalized volume}

We begin with a geometric lemma. We note that a {\em geodesic metric space} is a metric space $(X,d)$ where the distance between two points is attained by the length of a path between the points.

\begin{lemma}\label{path_fraction}
Let $Z$ be a collection of points in  a geodesic metric space $(X,d)$ such that for any collection of  $N+1$ points in $Y$ there are two that are at least a distance $\delta$ apart. Let
$$\alpha\colon [0,1]\to X$$
be a rectifiable path and let $L_\epsilon(\alpha)$ be the length of the path that is disjoint from the $\epsilon$-neighborhood of $Z$. Then for $\epsilon < \frac{\delta}{2N}$ 
$$L_\epsilon(\alpha) \ge \frac{\delta-2N\epsilon}{\delta}\left( d(\alpha(0),\alpha(1))-2N\epsilon\right).$$
\end{lemma}

{\bf Proof:} 
For each $z \in Z$ let
$$U_z = \{t \in [0,1]\ |\ d(\alpha(t), z) < \epsilon\}$$
and let $U$ be the union of the $U_z$. Note that for any $t\in [0,1]$ there are at most $N$ points $z \in Z$ such that $\mathcal N_{\epsilon}(z)$ intersects  the $(\delta-2\epsilon)/2$-neighborhood of $\alpha(t)$ and therefore there is a neighborhood of $t$ that intersects at most $N$ of the $U_z$. As $[0,1]$ is compact this implies that there are finitely many $z \in Z$ with $U_z\neq \emptyset$.

We claim we that we can find $z_1, \dots, z_m$ in $Z$ and 
$$0=t^+_0\le t^-_1 < t^+_1\le t^-_2 < \dots \le t^-_m < t^+_m \le 1=t^-_{m+1}$$
such that 
\begin{itemize}
\item $t^-_i \in \bar U_{z_i}$;
\item $t^+_i = \sup U_{z_i}$;
\item $\alpha([t^+_{i-1},t^-_i])$ is disjoint from $\mathcal N_\epsilon(Z)$.
\end{itemize}
We assume that the first $i$ points and values have been chosen and then find $z_{i+1}$ and $t^{\pm}_{i+1}$. Let $t^-_{i+1}$ be the infimum of $(t^+_i, 1]\cap U$. As there are finitely many non-empty $U_z$, there must be some $z \in Z$ with $t^-_{i+1}$ the infimum of $(t^+_i, 1]\cap U_z$. We let $z_{i+1} = z$ and $t^+_{i+1} = \sup U_z$. This process terminates (and $m=i$) when either $(t^+_i, 1]\cap U$ or $t^+_i = 1$.

Note that this implies that $d(\alpha(t^-_i), \alpha(t^+_i)) \le 2\epsilon$ and
$$\sum d(\alpha(t^+_{i-1}), \alpha(t^-_i)) \le L_\epsilon(\alpha).$$
Therefore
\begin{eqnarray*}
d(\alpha(0), \alpha(1)) & \le & \sum d\left(\alpha(t^-_i), \alpha(t^+_i\right)) + \sum d\left(\alpha(t^+_{i-1}), \alpha(t^-_i\right))\\
& \le & 2m\epsilon + L_\epsilon(\alpha).
\end{eqnarray*}

We need to show that $2m\epsilon$ is only a controlled portion of $d_(\alpha(0), \alpha(1))$. For this we choose a non-negative integer $k$ such that $kN < m \le (k+1)N$. 
Then we let $j_1$ be the smallest index such that there exists an $i_1< j_1$ with $d(z_{i_1}, z_{j_1}) \ge {\delta}$. Note that $j_1 \le N+1$. Then, as above,
$$\delta -2\epsilon \le d(\alpha(t^+_{i_1}), \alpha(t^-_{j_1})) \le 2(N-1)\epsilon + L_\epsilon\left(\alpha|_{[t^+_{i_1}, t^-_{j_1}]}\right).$$
Repeating this argument we get $i_\ell$ and $j_\ell$ for $\ell = 1, \dots, k$ where $j_{\ell-1} \le i_\ell < j_\ell$, $j_\ell - j_{\ell-1} \le N$ and
$$\delta -2\epsilon \le d(\alpha(t^+_{i_\ell}), \alpha(t^-_{j_\ell})) \le 2(N-1)\epsilon + L_\epsilon\left(\alpha|_{[t^+_{i_\ell}, t^-_{j_\ell}]}\right).$$
Summing these inequalities and rearranging we get
$$k \le \frac{L_\epsilon(\alpha)}{\delta-2N\epsilon}.$$
As $m \le (k+1)N$ our previous bound on $d(\alpha(0), \alpha(1))$ becomes
$$d(\alpha(0), \alpha(1)) \le 2(k+1)N\epsilon + L_\epsilon(\alpha).$$
Combining the two inequalities and rearranging gives the result.
\eproof

\begin{lemma}\label{flow_bound}
Assume that $0<\epsilon \le \epsilon_0$ and let $Y_t$ be a path on $\overline{\Teich(S)}$ such that   on $E=\{ t\ | \ d_{\rm WP}(Y_t, \mathcal G) > \epsilon\}$ the path is smooth and the tangent vector is the Weil-Petersson gradient of $-V_R$ and for $[u,v]$  a connected component of the path $Y_t$ in $E^c$ we have $V_R(Y_v) \leq V_R(Y_u)$. Then
$$V_R(Y_u) - V_R(Y_v) \ge 
A(\epsilon, S)\frac{\delta_0-2^{n(S)+1}\epsilon}{\delta_0} \left(d_{\rm WP}(Y_a,Y_b) - 2^{n(S)+1}\epsilon\right).$$

\end{lemma}

{\bf Proof:} Let
We have that $E$ is a collection $\mathcal I$ of open intervals. By assumption, for $t\in E$ the tangent vector $\dot Y_t$ of $Y_t$ is the Weil-Petersson gradient of $-V_R$ so by Theorem \ref{variational},
$$\|\dot Y_t\|_{\rm WP} = \|\phi_{Y_t}\|_2.$$
By Theorem \ref{nearnode} we also have that for $t\in E$,
$$\|\phi_{Y_t}\|_2 \ge A(\epsilon, S).$$
Again applying the variational formula, Theorem \ref{variational}, to an interval $(s,t)$ in $\mathcal I$ we have
\begin{eqnarray*}
V_R(Y_s) -V_R(Y_t) &=& \int_s^t \|\phi_{Y_t}\|_2^2 dt \\
& \ge& \int_s^t A(\epsilon, S)\|\phi_{Y_t}\|_2 dt \\
& = & A(\epsilon, S) L\left(Y_{(s,t)}\right)
\end{eqnarray*}
where $L(Y_{(s,t)})$ is the length of the path from $s$ to $t$.
For any interval $[u,v]$ in $E^c$, by assumption we have $V_R(Y_u)-V_R(Y_v) > 0$. Therefore we have
$$ V_R(Y_a) - V_R(Y_b) \ge \sum_{(s,t) \in \mathcal I} V_R(Y_s) - V_R(S_t)$$
and therefore 
$$V_R(Y_a) - V_R(Y_b) \ge A(\epsilon, S) L(Y_{\mathcal I})$$
where 
$$L(Y_{\mathcal I}) = \sum_{(s,t) \in \mathcal I}  L\left(Y_{(s,t)}\right).$$

For any collection of $2^{n(S)}+1$ simplices in $\mathcal C(S)$ there must be at least two that contain intersecting curves. Therefore by Theorem \ref{node_bound} for any collection of $2^{n(S)}+1$ points in $\mathcal G=\mathcal G(N;S,X)$ there are at least two that are a distance $\delta_0$ apart in the Weil-Petersson metric on $\overline{\Teich(S)}$ and we can apply Lemma \ref{path_fraction} with $Z = \mathcal G$  the set of points and $N = 2^{n(s)}$. Noting that $L_\epsilon(Y_{[a,b]})  = L(Y_{\mathcal I})$ by Lemma \ref{path_fraction} we have
$$L(Y_{\mathcal I}) \ge \frac{\delta_0 - 2^{n(S)+1}\epsilon}{\delta_0} \left(d_{\rm WP}(Y_a,Y_b) - 2^{n(S)+1}\epsilon\right).$$
Combining this with our above bound on the differences between renormalized volumes gives the result. \eproof

\subsubsection*{Convergence in the Weil-Petersson Completion}

\begin{prop}\label{flow_converges}
Let $Y_t$ be a flow line of the Weil-Petersson gradient flow of $-V_R$. Then $Y_t$ converges in $\overline{\Teich(S)}$ to a $\hat Y \in \mathcal G$.
\end{prop}

{\bf Proof:} By Lemma \ref{flow_bound} for every positive distance $d>0$ there is a $v>0$ such that if $d_{\rm WP}(Y_s, Y_t) \ge d$ then $V_R(Y_s) - V_R(Y_t) \ge v$. Renormalized volume is bounded below (and is in fact non-negative) and therefore $V_R(Y_t)$ converges as $t\rightarrow \infty$. In particular there exists a $T>0$ such that if $s,t >T$ then $V_R(Y_s) - V_R(Y_t) < v$ and $d_{\rm WP}(Y_s, Y_t) < d$. It follows that $Y_t$ converges in $\overline{\Teich(S)}$ as $t\rightarrow\infty$.

The lower bound on renormalized volume also implies that the integral
$$\int_0^\infty \|\phi_{Y_t}\|_2^2 dt < \infty.$$
Therefore we can find a sequence $t_i$ with $\|\phi_{Y_{t_i}}\|_2\to 0$ as $i\to\infty$. Theorem \ref{nearnode} then implies that any accumulation point of the sequence will lie in $\mathcal G$. As we have just seen that the entire path converges, this implies that the limit of $Y_t$ as $t\to \infty$ lies in $\mathcal G$.
\eproof

\subsubsection*{The surgered flow}\label{surgery-section}

\begin{prop}
Fix $\epsilon>0$. For all $Y \in \Teich(S)$ there exists a path $Y_t$ in $\overline{\Teich(S)}$ with $Y= Y_0$ such that 
\begin{itemize}
\item On $\{t \ | \ d_{\rm WP}(Y_t, \mathcal G) > \epsilon\}$ the path is smooth and the tangent vector is the Weil-Petersson gradient of $-V_R$;
\item If $a < b$  and $[a,b]$ is a connected component of $\{t \ |\ d_{\rm WP}(Y_t, \mathcal G) \leq \epsilon\}$ then $V_R(Y_b) < V_R(Y_a)$.
\item $Y_t \to Y_\geod$ as $t\to \infty$.
\end{itemize}
\end{prop}

{\bf Proof:} We claim there exists an integer $k\ge 0$ such that for $i=0, \dots, k$ there are a family of paths $Y^i_t$ and simplices $\tau_0, \dots, \tau_k$ in $\mathcal C(S)$ with
\begin{itemize}
\item $Y= Y^i_0$;

\item $Y^i_t$ passes through $Y_{\geod}^{\tau_0}, \dots, Y_{\geod}^{\tau_{i-1}}$;

\item $V_R(Y_{\geod}^{\tau_{j-1}}) < V_R(Y_{\geod}^{\tau_{j}})$ for $j=1,\ldots i-1$;

\item if $d_{\rm WP}(Y^i_t, \mathcal G) > \epsilon$ the path is smooth and the tangent vector $\dot Y^i_t$ is the Weil-Petersson gradient of $-V_R$;

\item $Y^i_t \to Y^{\tau_i}_\geod$ as $t\to \infty$ and $\tau_k = \emptyset$.
\end{itemize}
We start by letting $Y^0_t$ be the flow line of the Weil-Petersson gradient of $-V_R$ with $Y^0_0 = Y$. By Proposition \ref{flow_converges}, there is a simplex $\tau_0$ in $\mathcal C(S)$ such that $Y_t$ converges to some $Y^{\tau_0}_\geod\in \mathcal G$ where $\tau_0$ are the nodes of $Y^{\tau_0}_\geod$. 

Now assume $Y^0_t, \dots, Y^i_t$ and $\tau_0, \dots, \tau_i$ have been chosen. If $\tau_i =\emptyset$ then $k=i$ and we are done. If not, we form $Y^{i+1}_t$ as follows.
As $Y^i_t\to \tau_i$ there exists a $t_0$ such that if $t>t_0$ then $d_{\rm WP}(Y_t, Y^{\tau_i}_\geod) < \epsilon/2$. By Proposition \ref{no_local_min}, there is a path $Z_t$ with $Z_0 = Y^{\tau_i}_\geod$, $Z_t \in \Teich(S)$ and $V_R(Z_t) < V_R(Y^{\tau_i}_\geod)$. We can then choose $t_1$ such that if $0<t<t_1$ then $d_{\rm WP}(Y^{\tau_0}_\geod, Z_t) < \epsilon/2$. We then define $Y^1_t$ by
\begin{itemize}
\item $Y^{i+1}_t = Y^i_t$ if $t \le t_0$;

\item $Y^{i+1}_{[t_0, t_0+1)}$ is a reparameterization of $Y^i_{[t_0, \infty)}$;

\item $Y^{i+1}_t = Z_{t-t_0-1}$ if  $t \in [t_0+1, t_0 + t_1 + 1]$;

\item for $t\ge t_0 + t_1 +1$, $Y^{i+1}_t$ is a flow line of the Weil-Petersson gradient of $-V_R$.
\end{itemize}
For large $t$, $Y^{i+1}_t$ is a gradient flow line so once again by Proposition \ref{flow_converges}, $Y^{i+1}_t \to Y^{\tau_{i+1}}_\geod\in \mathcal G$ where curves in the simplex $\tau_{i+1}$ are the nodes of $Y^{\tau_{i+1}}_\geod$.

We now show that the process terminates. Observe that $V_R(Y^{\tau_{i+1}}_\geod) < V_R(Y^{\tau_i}_\geod)$ as the path $Y^{i+1}_t$ passes through $Y^{\tau_i}_\geod$, $V_R(Y^{i+1}_t)$ is decreasing, and $V_R(Y^{i+1}_t) \to V_R(Y^{\tau_{i+1}}_\geod)$ as $t\to\infty$ by Corollary \ref{extension}. This implies that all of the $\tau_i$ are distinct and $V_R(Y^{\tau_i}_\geod)$ is decreasing in $i$.

The flows $Y^i_t$ satisfy the conditions of Theorem \ref{flow_bound} so there exists a $v= v(\epsilon, \delta_0) >0$ such that if $d_{\rm WP}(Y^i_a, Y^i_b) \ge \delta_0$ then $V_R(Y^i_a) - V_R(Y^i_b) \ge v$. As we noted above, for any collection of $2^{n(S)}+1$ simplices in $\mathcal C(S)$ there will be at least two that contain intersecting curves. Therefore for any $i\ge 0$ there exist $j< \ell$ in $\{i, \dots, i + 2^{n(S)}\}$ such that $\tau_j$ and $\tau_\ell$ contain intersecting curves. By Theorem \ref{node_bound}  we then have $d_{\rm WP}(Y^{\tau_j}_\geod, Y^{\tau_\ell}_\geod) \ge \delta_0$. As $Y^{i+2^{n(S)}+1}_t$ passes through  $\tau_j$ and $\tau_\ell$, in that order (with possibly $i=j$ or $\ell=i+2^{n(S)}$), we have
$$V_R\left(Y^{\tau_i}_\geod\right) - V_R\left(Y^{\tau_{i+2^{n(S)}}}_{\rm geod}\right) \ge V_R(Y^{\tau_j}_\geod) - V_R(Y^{\tau_\ell}_\geod) \ge v.$$
Therefore, if the paths are defined up to $i$ with $2^{n(S)}m \leq i \leq 2^{n(S)}(m+1)$ we have
$$V_R(Y) - V_R(Y^{\tau_i}_\geod) \ge V_R(Y^{\tau_0}_\geod) - V_R(Y^{\tau_i}_\geod) \ge mv.$$
As $V_R \ge 0$ this implies that
$$i \le 2^{n(S)}\left(\frac{V_R(Y)}{v}+1\right).$$
Therefore the process must terminate.
\eproof.

We now use the above to give a new proof of the following theorem of Storm.

\begin{cor}[Storm, {\cite{storm2,storm1}}]\label{cvol_inf}
Let $N$ be a compact hyperbolizable acylindrical 3-manifold without torus boundary components. Then $V_C$  has a unique minimum at the structure $M_{\rm geod} \in CC(N)$ with totally geodesic convex core boundary.
 \end{cor}

The minimality of $M_{\rm geod}$ was the main result in \cite{storm2} and the uniqueness is a corollary of the main result in \cite{storm1} which considers the general case of $N$ with incompressible boundary.

{\bf Proof:} 
Let $Y \neq Y_{\rm geod}$, then using surgered flow,  we have the path  $Y_t$ with $Y_t \in \overline{\Teich(\partial N)}$ from $Y$ to $Y_{\rm geod}$ with $V_R(M_Y) > V_R(M_\geod)$.  Therefore
$$V_C(M_Y) \geq V_R(M_Y) > V_R(M_{\rm geod}) = V_C(M_{\rm geod}).$$
Thus $V_C$ has unique minimum at $M_{\rm geod}$.
\eproof

Note that in the course of the proof we have shown that the unique minimum of $V_R$ also occurs at $M_\geod$. In the relatively acylindrical case, we no longer have $V_C(M_\geod) = V_R(M_\geod)$ but otherwise the above proof goes through to give the following more general version of Storm's theorem for renormalized volume.

\begin{cor}{}\label{rvol_inf}
Let $(N;S)$ be a compact hyperbolizable relatively acylindrical 3-manifold without torus boundary components. Then $V_R$  has a unique minimum at the structure $M_{\rm geod} \in CC(N;S,X)$ with totally geodesic convex core boundary facing $S$.
 \end{cor}
 
 In \cite{BBB} we proved that Corollaries \ref{cvol_inf} and \ref{rvol_inf} are equivalent. Here we are directly proving both statements. A version of Corollary \ref{rvol_inf} was also proved by Pallete (\cite{PAL}) using different methods.

Also applying Lemma \ref{flow_bound} to the surgered flow path gives
\begin{theorem}
For all $\epsilon\le \epsilon_0$
$$V_R(Y) - V_R(Y_{\rm geod}) \ge A(\epsilon, S) \frac{\delta_0 - 2^{n(S)+1}\epsilon}{\delta_0} \left(d_{\rm WP}(Y, Y_{\rm geod})- 2^{n(S)+1}\epsilon\right).$$
\end{theorem}

Theorem \ref{main bound} then follows  from the above by choosing $\epsilon = \min(\epsilon_0, \delta_0/2^{n(s)+2})$ and letting 
$$A(S) = \frac{1}{2}A(\epsilon, S) \qquad\mbox{and} \qquad \delta = \frac{\delta_0}{2}.$$

\newcommand{\U}{{\mathbb U}}

We also recall Schlenker's upper bounds. His argument was originally for quasifuchsian manifolds but as we will see it holds whenever $(N;S)$ has incompressible boundary.
\begin{theorem}\label{upper bounds} Let $(N;S)$ have incompressible boundary. Then
$$|V_R(Y) - V_R(Y')| \le 3\sqrt{\frac{\pi}{2}|\chi(\partial N)|}d_{\rm WP}(Y, Y')$$
\end{theorem}

{\bf Proof:}
As noted in Corollary \ref{extension} the norm of the Weil-Petersson gradient of $V_R$ is bounded above by $\frac32\sqrt{\area(Y)} = 3\sqrt{\frac\pi 2 |\chi(S)|}$. Integrating this bound along a Weil-Petersson geodesic segment from $Y$ to $Y'$ gives the result. \eproof

We can now use the above  to prove Theorem \ref{qfmain} which we now restate.\newline

\noindent{\bf Theorem \ref{qfmain}}
{\em Let $S$ be a closed surface of genus $g\ge 2$. Then
$$A(S)(d_{\rm WP}(X,Y) - \delta) \le V_C(Q(X,Y)) \le 3\sqrt{\frac{\pi}2 |\chi(S)|}d_{\rm WP}(X,Y) + 6\pi|\chi(S)|.$$
}

{\bf Proof:}
If $N = S\times [0,1]$ then a {\em Bers slice} is the deformation space $CC(N; S\times\{0\},X)$ where $X$ is a fixed conformal structure on $S$. Manifolds in this deformation space are {\em quasifuchsian} and the manifold $M_Y \in CC(N; S\times\{0\},X)$ in our general notation is usually referred to as $Q(X,Y)$.

We apply  Theorem \ref{main bound} to this case.  Then $Q(X,X)$ is the Fuchsian manifold so $Y_\geod = X$ and $V_R(Y_\geod) = 0$.   Therefore we have
$$ 
A(S)(d_{\rm WP}(X,Y) -\delta)
\le \Rvol(Q(X,Y)) .$$

Combining this lower bound with the bound of  Schlenker (see  \cite[Theorem 1.2]{schlenker-qfvolume}), we have 
$$ 
A(S)(d_{\rm WP}(X,Y) -\delta)
\le \Rvol(Q(X,Y)) \le 
3\sqrt{\frac{\pi}{2}|\chi(S)|} d_{\rm
  WP}(X,Y) 
$$
By \cite{BBB}, for any convex cocompact $M$ 
$$V_R(M) +\frac{1}{4}L(\beta_M) \leq V_C(M) \leq V_R(M) +\frac{1}{2}L(\beta_M).$$
Also  for $\partial N$ incompressible $L(\beta_M) \leq 6\pi |\chi(\partial N)|$ (see \cite{BBB}).
Therefore the result follows.
\eproof

Theorem \ref{maincorevolume}  follows identically as in the  proof of Theorem \ref{qfmain} above.

%% file: appendix.tex
%% !TEX root =renorm_wp.tex

\appendix
\section{Appendix: A Weil-Petersson estimate}

We recall that the Margulis constant in two dimensions is $\epsilon_2 = \sinh^{-1}(1)$. In this section we prove the following proposition:
\begin{prop}\label{wp estimate}
Let $\tau$ be a simplex in ${\mathcal C}(S)$ and $Y\in \Teich(S)$ a hyperbolic surface such that for each curve $\beta\in \tau$, $\ell_\beta(Y) \leq \ell_0$ where $0 < \ell_0 < 2\epsilon_2$. Let $\hat Y \in \Teich(S_\tau)$ be such that the cover $\check Z$ of $Y$ associated to $S\backslash \tau$ conformally embeds in $\hat Y$. Then
$$d_{\rm WP}(Y, \hat Y)\le 2\pi\sqrt{\frac{2\sinh(\ell_0/2)}{\ell_0(1-\sinh(\ell_0/2))}} \sqrt{\ell_\tau(Y)}.$$
\end{prop}

We will use the following criteria for convergence in the Weil-Petersson completion: Let $\tau$ be a simplex in $\mathcal C(S)$ and $\hat Y$ a surface in $\Teich(S_\tau)$. Then a sequence $Y_i\in\Teich(S)$ converges to $\hat Y$ in $\overline{\Teich(S)}$ if for all simple closed curves $\gamma$ with $i(\gamma, \tau) = 0$ we have $\ell_\gamma(Y_i) \to \ell_\gamma(\hat Y)$. In particular the length of the curves in $\tau$ must converge to zero. We will use the following lemma to verify this criteria.
\begin{lemma}\label{WP limit}
Let $R\subset S$ be a proper, essential, non-annular subsurface of a finite type surface $S$. Let $R_i$ and $S_i$ be conformal structures on $R$ and $S$, respectively, such that there is a conformal embedding $R_i \hookrightarrow S_i$ in the homotopy class of $R\hookrightarrow S$. If $\ell_{\partial R}(R_i) \to 0$ then for all simple closed curves $\gamma$ on $R$ we have
$$\underset{i\to\infty}{\lim} \ell_\gamma(R_i) = \underset{i\to\infty}{\lim} \ell_\gamma(S_i)$$
where the lengths are measured on the completed hyperbolic metrics on the respective conformal structures.
\end{lemma}

{\bf Proof:} Let $R^\gamma_i$ and $S^\gamma_i$ be the annular covers of $R_i$ and $S_i$ corresponding to the curve $\gamma$.  Then there is a conformal embedding $R^\gamma_i \hookrightarrow S^\gamma_i$ that is a homotopy equivalence. Therefore
$$ \frac{\pi}{\ell_\gamma(R_i)} = m\left(R^\gamma_i\right) \le m\left(S^\gamma_i\right) = \frac{\pi}{\ell_\gamma(S_i)}$$
where $m(\cdot)$  is the modulus of the annulus.

To get a bound in the other direction we let $D_i$ be the distance, in the $S_i$-metric, from the geodesic representative of $\gamma$ in $S_i$ to the complement of $R_i$ and denote the $D_i$-neighborhood of the geodesic core of $S^\gamma_i$  by $S^\gamma_i(D_i)$. Then $S^\gamma_i(D_i)$ will be contained $R_i$ and it follows that
$$m\left(S^\gamma_i(D_i)\right)= \frac{\pi - \epsilon_i}{\ell_{\gamma}\left(S_i\right)} \le m(R^\gamma_i)$$
where $\epsilon_i$ only depends on $D_i$ and $\epsilon_i\to 0$ as $D_i\to \infty$. To finish the proof we need to show that $D_i\to\infty$.

Let $C(R_i)$ be the convex core of $R_i$ and assume that each component of the boundary of $C(R_i)$ has length $< 2\epsilon_2$. Then each component of the boundary of $C(R_i)$ will lie in the standard collar of the associated geodesic in $S_i$. As the length of the boundary curves of $C(R_i)$ limits to zero the depth of these curves in the standard $S_i$-collars will limit to infinity. In particular, the distance of any point in the $R$-component of the complement of the $S_i$-collars from the complement of the $R_i$ will also limit to infinity. As the geodesic representative of $\gamma$ in $S_i$ will be in this complementary region we have that $D_i\to\infty$, as desired.
\eproof

Let $A$ be a conformal annulus with finite modulus $m(A)$. Then $A$ can be realized as the quotient of the strip
$$S = \{ z\in \mathbb C | 0 < \Im s < \pi\}$$
by the translation
$$z \mapsto z+ \frac \pi {m(A)}.$$
Define Beltrami differentials $\mu^t_A$ and $\mu^h_A$ so that their lifts to $S$ are $\tilde{\mu}^t_A = 1$ and $\tilde\mu^h_A = \sin^2 y$, respectively. Then $\mu$ is a {\em Teichm\"uller differential} on $A$ if it is a constant multiple of $\mu_A^t$ and is a {\em harmonic differential} on $A$ if it is a constant multiple of $\mu^h_A$.

\begin{lemma}\label{harmonic}
Let $\mu$ be a Beltrami differential on $Y$ such that on an annulus $A$, $\mu = c \mu^t_A$ is a Teichm\"uller differential. Assume that $\nu$ is the Beltrami differential with $\nu = 2c \mu^h_A$ on $A$ and $\nu = \mu$ on the complement of $A$. Then $\mu - \nu$ is an infinitesimally trivial Beltrami differential.
\end{lemma}

{\bf Proof:} We need to show that for any holomorphic quadratic differential $\phi \in Q(Y)$ the pairing of $\phi$ with $\mu-\nu$ is zero. The difference $\mu-\nu$ is supported on $A$ so our computation will be on fundamental domain in $S$ for the action $z\mapsto z+\pi/m(A)$. 
The restriction of $\phi$ to $A$ lifts to a holomorphic quadratic differential $g(z) dz^2$ on $S$ where $g$ is a periodic holomorphic function. That is
$$g(z + \pi/m(A)) = g(z).$$

Let
$$b(y) = \int_0^{\pi/m(A)} g(x+iy) dx.$$
If $Q$ is a rectangle whose top and bottom sides are horizontal segments from $x=0$ to $x= \pi/m(A)$ at heights $y_0<y_1$ then
$$\int_{\partial Q} g(z) dz = b(y_0) - b(y_1)$$
since the periodicity of $g(z)$ implies that the line integrals over the vertical sides cancel. As $g(z)$ is holomorphic the line integral around $\partial Q$ is zero and therefore $b(y_0) = b(y_1)$ which implies that $b(y) \equiv b$ is a constant function.

Using this we now compute the pairing:
\begin{eqnarray*}
\int_Y (\mu-\nu) \phi &=& \int_A (\mu-\nu) \phi\\ &=& \int_0^\pi\int_0^{\pi/m(A)}c(1-2\sin^2y) g(x+iy) dx dy\\
&=& \int_0^\pi cb (1-2\sin^2 y) dy \\
& =& 0.
\end{eqnarray*}
\eproof

In practice it is easier to construct deformations where the tangent vectors are infinitesimal Teichm\"uller differentials on annuli. We can use the previous lemma to bound the Weil-Petersson norm of these deformations.
\begin{lemma}\label{harmonic bound}
Let $A_i$ be a collection of disjoint annuli on $Y$ with finite moduli $m_i$. If
$$\mu = \sum_i c_i \mu^t_{A_i}$$
is a Beltrami differential on $Y$  then
$$\|[\mu]\|_2^2 \le 2\pi^2 \sum \frac{|c_i|^2}{m_i}.$$
\end{lemma}

{\bf Proof:} By Lemma \ref{harmonic} the Beltrami differential $\mu$ is equivalent to
$$\nu =  2\sum_i c_i \mu^h_{A_i}$$
so
$$\|[\mu]\|_2^2 = \|[\nu]\|_2^2  \le  \int_Y \|\nu\|^2 da_Y$$
where $da_Y$ is the area form for the hyperbolic metric on $Y$. By the Schwarz lemma if $da_i$ is the area form for the complete hyperbolic metric on $A_i$ then $da_Y < da_i$. On the strip $S$ the area form $da_i$ lifts to $\frac1{\sin^2y} dx dy$ so
\begin{eqnarray*}
\int_Y \|\nu\|^2 da_Y & \le & 4\sum_i \int_{A_i} |c_i|^2\|\mu^h_{A_i}\|^2 da_i\\
& = & 4\sum_i |c_i|^2 \int_0^\pi\int_0^{\pi/m(A_i)} \frac{(\sin^2 y)^2}{\sin^2 y} dx dy\\
& = & 4\sum \frac{|c_i|^2 \pi^2}{2m_i}.
\end{eqnarray*}
\eproof

We can now describe the strategy of the proof of Proposition \ref{wp estimate}. Let $Z\subset Y$ be the complement of the geodesic representatives of $\tau$ in $Y$. Then $Z$ will lift to $\check Z$  and conformally embed in both $Y$ and $\hat Y$. We will construct a family of quasiconformal deformations of $\hat Y$ to itself where the tangent vectors of these deformations will be Teichm\"uller differentials on a collection of annuli that lie in $Z\subset \hat Y$. As $Z$ is also a subsurface of $Y$ this will define a family of quasi-conformal deformations of $Y$, but here the surface will change along the deformation. This will define a path in $\Teich(S)$. We will use Lemma \ref{WP limit} to see that this path converges to $\hat Y$ and Lemma \ref{harmonic bound} to bound above the Weil-Petersson length of the path.

\subsubsection*{The cusp deformation}
Every cusp $\frak C$ of a hyperbolic Riemann surface can be parameterized as the quotient of the horodisk
$$\frak H = \{ z \in \mathbb C | \Im z \ge 1\}$$
by the translation 
$$z\mapsto z+2.$$
If we let
$$\frak H(m) = \{ z \in \mathbb C |  1 \le \Im z \le 2m+1\}$$
then the quotient $\frak C(m)$ of $\frak H(m)$ is an annulus of modulus $m$. Define maps
$$f^m_t\colon \frak C \to \frak C$$
such that $f^m_t$ is
\begin{itemize}
\item constant in the $x$-variable;

\item an affine map from $\frak C(m)$ to $\frak C(e^t m)$;

\item is conformal in the complement of $\frak C(m)$.
\end{itemize}
At time $t$ the infinitesimal Beltrami differential $\nu_t$ for this path will be supported on the annulus $\frak C(e^tm)$ and using the fact that 
$$f^{e^tm}_s \circ f^m_t = f^m_{s+t}$$
we see that the lift of $\nu_t$ to $\frak H$ is supported on $\frak H(e^tm)$ with  $\tilde\nu_t = -\frac12$. In particular, $\nu_t$ is a Teichm\"uller differential on $\frak C(e^t m)$.

\subsubsection*{The deformation of $\hat Y$ and $Y$}
Each curve in $\tau$ is a node of $\hat Y$ and there are two associated cusps in $\hat Y$. If $\tau$ has $k$ curves we label the two cusps associated to the $i$th node by $\frak C_i^\pm$ and assume that the modulus $m_i$ has been chosen such that the annuli $\frak C_i^\pm (m_i)$ lie in $Z$.

With this choice of moduli we define a family of maps
$$f_t \colon \hat Y \to \hat Y$$
by setting $f_t$ to be the map 
$f^{m_i}_t$ on the cusps $\frak C^\pm_i$ and to be the identity on the complement of the cusps. (There may be cusps of $\hat Y$ that don't correspond to nodes in $\tau$. The map is the identity here.)  The Beltrami differentials $\mu_t$ for this family of maps are supported on the annuli $\frak C^\pm_i(m_i)$. As these lie in $Z$, the $\mu_t$ are also a family of Beltrami differentials on $Y$ so we have two one-parameter families of surfaces $Z_t$ and $Y_t$ with $Z_t$ conformally embedding in $Y_t$. The $Z_t$ also conformally embed in $\hat Y$.

\medskip

{\bf Proof of Proposition \ref{wp estimate}:}
Let $\beta_i$ be the $i$th curve of $\tau$ and let $\beta^\pm_i$ be the two curves that are homotopically distinct in $S\backslash \tau$ but are both homotopic in $S$ to $\beta_i$. Let $Z_t^{\beta^\pm_i}$ be the annular cover of the component of $Z_t$ containing ${\beta^\pm_i}$. Then
$$m\left(Z^{\beta^\pm_i}_t\right) \ge e^tm_i$$
and therefore
$$\ell_{\partial Z_t}(Z_t) \to 0$$
as $t\to \infty$. By Lemma \ref{WP limit}  for all non-peripheral simple closed curves $\gamma$ in $R$ we  have
$$\underset{t\to\infty}{\lim} \ell_\gamma(Z_t) = \underset{t\to\infty}{\lim} \ell_\gamma(Y_t)$$
and
$$\underset{t\to\infty}{\lim} \ell_\gamma(Z_i) = \underset{t\to\infty}{\lim} \ell_\gamma(\hat Y) = \ell_\gamma(\hat Y).$$
It follows that
$$\underset{t\to\infty}{\lim} \ell_\gamma(Y_t) = \ell_\gamma(\hat Y)$$
so $Y_t\to \hat Y$ in $\overline{\Teich(S)}$.

The tangent vector of the path are Teichm\"uller differentials on $2k$ disjoint annuli with coefficients $-\frac12$. At time $t$, two of these annuli have modulus $e^t m_i$ so integrating the estimate from Lemma \ref{harmonic bound} we have
$$d_{\rm WP}(Y, \hat Y) \le \int_0^\infty \sqrt{\pi^2 \sum_i \frac1{m_i e^t}} = 2\pi \sqrt{\sum\frac1{m_i}}$$

To finish the proof we need to bound the $m_i$ from below. As $\check Z$ is a cover of $Y$, $\ell_{\beta_i^\pm}(\check Z)  = \ell_{\beta_i}(Y)$. By the Schwarz lemma, the geodesic representative of $\beta_i^\pm$ in $\check Z$ will lie in the $\ell_{\beta_i^\pm}(\check Z)/2$-thin part of the associated cusps $\frak C^\pm$ of $\hat Y$. If $p \in \frak C$ is a point in  our standard model  of a cusp with pre-image $z=x+iy\in \frak H$ then injectivity radius satisfies the formula
$$\sinh(\inj(p)) = \frac1y.$$
Note that while $z$ is not uniquely determined, the $y$-coordinate is. This implies that $\check Z$ will contain the annuli $\frak C(m_i)$ where
$$m_i =  \frac{1}{2}\left(\frac{1}{\sinh(\ell_{\beta_i}(Y)/2)} - 1\right).$$
With our assumption that $\ell_{\beta_i^\pm}(\check Z) = \ell_{\beta_i}(Y) \le \ell_0$ we have
$$\sinh(\ell_{\beta_i}(Y)/2)  \le \frac{\sinh(\ell_0/2)}{\ell_0/2}\cdot \frac{\ell_{\beta_i}(Y)}2$$
and therefore
\begin{eqnarray*}
m_i &\ge& \frac{\ell_0}{2\sinh(\ell_0/2)\ell_{\beta_i}(Y)} - \frac12 \\
& = & \frac{\ell_0- \sinh(\ell_0/2)\ell_{\beta_i}(Y)}{2\sinh(\ell_0/2)\ell_{\beta_i}(Y)}\\
&\ge &\frac{\ell_0- \sinh(\ell_0/2)\ell_0}{2\sinh(\ell_0/2)\ell_{\beta_i}(Y)} =  \frac{\ell_0(1-\sinh(\ell_0/2))}{2\sinh(\ell_0/2)\ell_{\beta_i}(Y)} .
\end{eqnarray*}
It follows that
\begin{eqnarray*}
d_{\rm WP}(Y, \hat Y) &\le& 2\pi \sqrt{\sum_i \frac{2\sinh(\ell_0/2)\ell_{\beta_i}(Y)}{\ell_0(1-\sinh(\ell_0/2))}}\\
& = & 2\pi\sqrt{\frac{2\sinh(\ell_0/2)}{\ell_0(1-\sinh(\ell_0/2))}} \sqrt{\ell_\tau(Y)}.
\end{eqnarray*}
\eproof